%% file: twomicro.tex
\documentclass{amsart}
\usepackage{standard}
\usepackage{epic,eepic}
\usepackage{verbatim}

\numberwithin{theorem}{section}

\newcommand{\clcon}{\mathcal{A}_{\CL}} 

\renewcommand{\Im}{\operatorname{Im}}
\DeclareMathOperator{\Opl}{{{}^{\mathit{h}}Op_l}}
\DeclareMathOperator{\Opr}{{{}^{\mathit{h}}Op_r}}
\DeclareMathOperator{\Opw}{{{}^{\mathit{h}}Op_W}}
\newcommand{\restrictedto}{\upharpoonright}
\newcommand{\FT}{\mathcal{F}} 
\newcommand{\tot}{\text{tot}}
\newcommand{\hI}{\hat{I}}

\newcommand\sphere{\mathbb{S}}

\newcommand{\TT}{\mathbb{T}}

\newcommand{\bdf}{\rho}
\newcommand{\sbdf}{\tilde\rho_{\ff}}

\newcommand{\U}{\mathcal{U}}
\newcommand{\cS}{\mathcal{S}}

\newcommand{\residual}{{\mathcal{R}}}
\DeclareMathOperator{\esssupp}{{ess supp}}
\newcommand{\symbsp}{\mathsf{S}}
\newcommand{\pssp}{\mathsf{S}_{0}} 
\newcommand{\sH}{\mathsf{H}}

\DeclareMathOperator{\liptic}{{ell}}
\newcommand{\ombar}{\overline{\omega}}
\newcommand{\CL}{{\text{cl}}}
\newcommand{\diag}{{\text{diag}}}
\newcommand{\loc}{{\text{loc}}}

\newcommand{\lag}{\mathcal{L}}
\newcommand{\Psitwo}{\Psi_{2,h}}
\newcommand{\Psih}{\Psi_h}
\newcommand{\Psism}{\tilde\Psi_h}

\newcommand{\ff}{{\text{ff}}}

\newcommand{\sidef}{{\text{sf}}}

\newcommand{\fcal}{\mathcal{F}}

\newcommand{\mcal}{\mathcal{M}}

\newcommand{\fcalh}{{}^h\mathcal{F}}

\DeclareMathOperator{\Id}{Id}

\DeclareMathOperator{\sub}{sub}

\DeclareMathOperator{\sigmatwo}{{}^2\sigma}
\DeclareMathOperator{\Sigmatwo}{{}^2\Sigma}
\DeclareMathOperator{\WFtwo}{{}^2WF}

\newcommand\Vf{{\mathcal V}}

\author{Andr\'as Vasy}
\address{Department of Mathematics\\Stanford University\\Stanford CA}
\email{andras@math.stanford.edu}
\author{Jared Wunsch}
\address{Department of Mathematics\\Northwestern University\\Evanston IL}
\email{jwunsch@math.northwestern.edu}
\thanks{The authors gratefully acknowledge financial support for this
project from the National Science Foundation, the first under
grant DMS-0201092, and the second under grant DMS-0700318.}

\title{Semiclassical second microlocal propagation of regularity and
integrable systems}

\begin{document}
\begin{abstract}
We develop a second-microlocal calculus of pseudodifferential operators in
the semiclassical setting.  These operators test for Lagrangian regularity
of semiclassical families of distributions on a manifold $X$ with respect
to a Lagrangian submanifold of $T^*X.$ The construction of the calculus,
closely analogous to one performed by Bony in the setting of homogeneous
Lagrangians, proceeds via the consideration of a model case, that of the
zero section of $T^*\mathbb{R}^n,$ and conjugation by appropriate Fourier
integral operators.  We prove a propagation theorem for the associated
wavefront set analogous to H\"ormander's theorem for operators of real
principal type.

As an application, we consider the propagation of Lagrangian regularity on
invariant tori for quasimodes (e.g.\ eigenfunctions) of an operator with
completely integrable classical hamiltonian.  We prove a secondary
propagation result for second wavefront set which implies that even in the
(extreme) case of Lagrangian tori with all frequencies rational, provided a
nondegeneracy assumption holds, Lagrangian regularity either spreads to
fill out a whole torus or holds nowhere locally on it.
\end{abstract}

\maketitle
\section{Introduction}
\subsection{Second microlocalization on a Lagrangian}
One purpose of the calculus of pseudodifferential operators is to test
distributions for regularity.  In the case of the semiclassical calculus,
regularity is measured by powers of the semiclassical parameter $h;$ if
$u_h$ is a family of distributions as $h\downarrow 0,$ one can, following
\cite{GS}, define a ``frequency set'' or (as we will refer to it here)
``semiclassical wavefront set'' inside the cotangent bundle $T^*X$ of the
underlying manifold $X,$ by decreeing that $p\notin \WF_h (u_h)$ if and only
if for arbitrary $k$ and $A_1,\dots,A_k \in \Psih^{-\infty}(X),$ microsupported
sufficiently close to $p,$ we have $h^{-k} A_1\dots A_k u_h \in L^2,$
uniformly as $h\downarrow 0$ (this uniformity will henceforth be tacit).
Here $\Psih^{-\infty}(X)$ stands for the algebra of smoothing semiclassical
pseudodifferential operators, of order $0$ in $h$ (thus uniformly bounded
on $L^2$); see Section~\ref{sec:calculus} for the details of the notation.
This ``oscillatory testing'' definition is quite flexible, and illustrates
the role of semiclassical pseudodifferential operators as test operators
for regularity relative to $L^2.$ With $\WF_h(u_h)$ also defined for
points at ``fiber infinity'' on the cotangent bundle, i.e.\ on
$S^*X=(T^*X\setminus o)/\RR^+$,
we have $\WF_h (u_h)=\emptyset$ if and
only if $u_h \in h^\infty L^2_{\loc}.$

Many distributions arising in the theory of PDE are, of course, not
$O(h^\infty)$ (or, in the conventional, homogeneous, theory, not
smooth); a great many of the examples that arise in practice, however, turn
out to be regular in a different way: they are \emph{Lagrangian
distributions,} associated to a Lagrangian submanifold $\lag \subset T^*X.$
We may characterize these distributions again by an ``iterated regularity''
criterion: if for all $k$ and all $A_1,\dots A_k \in \Psih^1(X),$ with
$\sigma_h(A_i)\equiv 0$ on $\lag$
$$
h^{-k}A_1\dots A_k u_h \in L^2,
$$ we say that $u$ is a \emph{Lagrangian distribution} with respect to
$\lag.$ This characterization, analogous to the Melrose-H\"ormander
characterization of ordinary (i.e.\ homogeneous, or non-semiclassical)
Lagrangian distributions, is equivalent to the statement that $u_h$ has an
oscillatory integral representation as a sum of terms of the form
$$
\int a(x,\theta,h) e^{i \phi(x,\theta)/h} \, d\theta,
$$ where $\phi$ parametrizes the Lagrangian $\lag$ appropriately (see, for
instance, \cite{Hormander:vol4} in the classical case, and
\cite{Alexandrova:Semiclassical} or \cite{Sjostrand-Zworski:Fractal} for an
account of semiclassical Lagrangian distributions).  We may, by limiting
the microsupport of the test operators $A_i,$ somewhat refine this
description of Lagrangian regularity to be \emph{local} on $\lag.$ It
remains, however, somewhat crude: it turns out to be quite natural to test
more finely, with semiclassical pseudodifferential operators whose
principal symbols are allowed to be singular at $\lag$ in such a way as to
be smooth on the manifold obtained by performing \emph{real blowup} on
$\lag$ inside $T^*X,$ i.e.\ by introducing polar coordinates about it.  The
resulting symbols localize not only on $\lag$ itself, but more finely, in
$SN(\lag),$ the spherical normal bundle.  (We may, by using the symplectic
structure, identify $SN(\lag)$ with $S(\lag),$ the unit sphere bundle of
$T(\lag),$ but we will not adopt this notation.)  The resulting
pseudodifferential calculus is said to be \emph{second microlocal}; there
is an associated wavefront set in $SN(\lag)$ whose absence (together with
absence of ordinary semiclassical wavefront set on $T^*X\backslash \lag$) is
equivalent to $u_h$ being a semiclassical Lagrangian distribution.  A
helpful analogy is that second-microlocal wavefront set in $SN(\lag)$ is to
failure of local Lagrangian regularity on $\lag$ as ordinary wavefront set
is to failure of local regularity on $X,$ better known as singular support.

The first part of this paper is devoted to the construction of the
semiclassical second microlocal calculus for a Lagrangian in $T^*X,$ and an
enumeration of its properties.  Other instances of second microlocalization
abound in the literature, although we know of none existing in the
semiclassical case, with respect to a Lagrangian.  Our approach stays
fairly close to that adopted by Bony \cite{Bony:Second} in the classical
case of homogeneous Lagrangians, and to that of Sj\"ostrand-Zworski
\cite{Sjostrand-Zworski:Fractal}, who construct a semiclassical second
microlocal calculus adapted to hypersurfaces in $T^*X.$

\subsection{An application to quasimodes of integrable Hamiltonians}

As an example of the power of second microlocal techniques in the
description of Lagrangian regularity, in the second part of the paper we
consider quasimodes of certain operators\footnote{In addition to the
condition that the principal symbol be real, we also require a hypothesis
on the subprincipal symbol: see \S\ref{section:propagation} for details.}
with real principal symbol with \emph{completely integrable} Hamilton flow.
Quasimodes are solutions to
$$
P_h u_h  \in h^k L^2
$$ for some $k$ (the order of the quasimode); eigenfunctions of
Schr\"odinger operators are of course motivating examples.  We further
assume that the foliation of the phsae space is (locally, at least) given
by \emph{compact} invariant tori; these tori are Lagrangian. (See
\cite{Arnold} for an account of the theory of integrable systems.)  The
Hamilton flow on a Lagrangian torus $\lag$ is given by quasi-periodic motion
with respect to a set of frequencies $\ombar_1,\dots, \ombar_n.$ It was
shown in \cite{Wunsch:Integrable} that local Lagrangian regularity on
$\lag$ propagates along Hamilton flow, hence if all frequencies are
irrationally related, it fills out the torus.  (The set on which local
Lagrangian regularity holds is open.)  Thus Lagrangian regularity is, on
one of these irrational tori, an ``all or nothing'' proposition: it obtains
either everywhere or nowhere on $\lag.$  In \cite{Wunsch:Integrable},
the opposite extreme case was also considered: Lagrangian tori on which
$\ombar_i/\ombar_j \in \QQ$ for each $i,j.$ Local Lagrangian regularity
must occur on unions of closed orbits, but in this case, these orbits need
not fill out the torus.  It was shown, however, that in the presence of a
standard nondegeneracy hypothesis (``isoenergetic nondegeneracy''), local
Lagrangian regularity propagates in one additional way: to fill in small
tubes of bicharacteristics.  This apparently mysterious and ungeometric
propagation phenomenon is elucidated here.  We study the propagation of
\emph{second microlocal} regularity on $SN(\lag),$ and find that it is invariant
under two separate flows: the Hamilton flow lifted to $SN(\lag)$ from the
blowdown map to $T^*X,$ and a second flow given by the next-order jets of
the Hamilton flow near $\lag\subset T^*X.$ This leads, in the case
considered in \cite{Wunsch:Integrable}, to the ``all or nothing'' condition
also holding for local Lagrangian regularity on nondegenerate rational
invariant tori: once again either the distribution is Lagrangian on the
whole of $\lag$ or nowhere locally Lagrangian on it.

The authors are very grateful to Maciej Zworski for helpful discussions
about this work, and in particular for directing them to reference
\cite{Hitrik-Sjostrand1}.

\section{The Calculus}\label{sec:calculus}
Let $X$ denote a manifold without boundary. We adopt the convention that
$\Psi_h^{m,k}(X)=h^{-k}\Psi_h^m(X)$ is the space of semiclassical
pseudodifferential operators on $X$ of differential order $m,$ hence given
locally by semiclassical quantization of symbols lying in $h^{-k}\CI([0,1)_h;
S^m(T^*X)).$ However, we almost exclusively
work microlocally in a compact subset of $T^*X\times 0\subset T^*X\times
[0,1)$, so the differential order, corresponding to the behavior of total
symbols at infinity in the fibers of the cotangent bundle, is irrelevant
for us, hence we also let
$$
\Psism^k(X)\subset h^{-k}\Psi_h^{-\infty}(X)=\Psi_h^{-\infty,k}(X).
$$ be the subalgebra consisting of ps.d.o's with total symbols compactly
supported in the fibers of $T^*X$ plus symbols in $h^\infty\CI([0,1);
S^{-\infty}(T^*X)).$ (For accounts of the semiclassical pseudodifferential
calculus, see, for instance, \cite{Martinez:book, Dimassi-Sjostrand,
Evans-Zworski}).  We will suppress the $h$-dependence of families of
operators (writing $P$ instead of $P_h$) of distributions (writing $u$
instead of $u_h$).

Let $\lag \subset T^*X$ be a Lagrangian submanifold with the restriction of
the bundle projection to $\lag$ being proper.  We will define a calculus of
pseudodifferential operators $\Psitwo(X;\lag)$ associated to $\lag$ with
the following properties.

\renewcommand{\theenumi}{\roman{enumi}}
\begin{enumerate}
\item
$\Psitwo^{*,*}(X;\lag)$ is a calculus: it is a bi-filtered algebra of
operators $A=A_h:\CI(X)\to\CI(X)$ with properly supported Schwartz
kernels,
closed under adjoints and asymptotic summation: if
$A_j\in\Psitwo^{m-j,l}(X;\lag)$ then there exists $A\in\Psitwo^{m,l}(X;\lag)$
such that $A-\sum_{j=0}^{N-1}A_j \in\Psitwo^{m-N,l}(X;\lag)$ for all $N$.

\item\label{principalsymbol} There is a principal symbol map
$$
\sigmatwo_{m,l}: \Psitwo^{m,l}(X;\lag) \to \clcon^{m-l}(\pssp),
$$ where $\pssp=[T^*X,\lag]$ denotes the \emph{real blowup} of $\lag$ as a
submanifold of $T^*X$ given by introducing normal coordinates about $\lag$
(see \cite{M} for extensive discussion or \cite{MR93d:58152} for a brief
account) and where
$$
\clcon^{m-l}(\pssp)
$$ is the space of classical conormal distributions with respect to
$SN(\lag),$ the spherical normal bundle of $\lag,$ which is canonically
identified with $\pa \pssp$: if {\em $\sbdf$ is a boundary defining
function for this face}, then\footnote{The requirement of compact support
which we have built into this definition is convenient but not strictly
necessary; see Remark~\ref{remark:noncompact}.}
$$
\clcon^{r}(\pssp)\equiv \sbdf^{r} \CI_c(\pssp).
$$
For brevity, we will let 
$$
S^r(\pssp) =  \clcon^{-r}(\pssp).
$$

The map $\sigmatwo$ is a $*$-algebra homomorphism, and fits into the short
exact sequence
\begin{equation}\label{ses}
0 \to \Psitwo^{m-1,l}(X;\lag) \to \Psitwo^{m,l}(X;\lag)
\overset{\sigmatwo_{m,l}}{\to} S^{l-m}(\pssp) \to 0.
\end{equation}
We remark in particular that the vanishing of the symbol only reduces the
order in one of the two indices.


\item
There is a quantization map
$$
\Op: \clcon^{-m,-l}(\symbsp)\to \Psitwo^{m,l}(X;\lag),
$$
where
$$
\symbsp=[T^*X\times [0,1); \lag\times 0]
$$ can be thought of as the space of ``total symbols'' of two-pseudors.
Note that $\pssp,$ the space on which principal symbols live, is one of the
boundary faces of the manifold with corners $\symbsp.$ Here again $\clcon$
refers to (compactly supported) classical conormal distributions, i.e.\
multiples of boundary defining functions times smooth functions on the
manifold with corners $\symbsp;$ the indices $-m,-l$ refer to the orders at
the front face of the blowup and the side face (i.e.\ the lift of $\pssp$)
respectively.

For brevity, we will let
$$
S^{m,l}(\symbsp) = \clcon^{-m,-l}(\symbsp).
$$

Since one boundary face of $\symbsp$ is $\pssp,$ if $a\in S^{l-m}(\pssp)$ we
may extend it to an element of $S^{0,l-m}(\pssp),$ and multiply by $h^{-m}$
to obtain
$\tilde a \in S^{m,l}(\symbsp).$ This we may quantize and obtain of course
$$
\sigmatwo_{m,l}(\Op(\tilde{a})) = a.
$$

\item
If $a \in S^{m,l}(\symbsp),$ let
  $\WF'(\Op(a))$ be defined as $\esssupp(a) \subset \pssp,$
where $\esssupp(a)^c$ is the set of points in
$\pssp\subset \symbsp$ which have a
neighborhood in which $a$ vanishes to infinite order {\em at}
$\pssp.$
 Then
  $\WF'$ in fact well-defined on $\Psitwo^{*,*}(X;\lag),$ and
$$
\WF'(A+B)\subset\WF'(A)\cup\WF'(B),\ \WF'(AB) \subset \WF'(A) \cap \WF'(B).
$$

For $A \in \Psitwo^{m,l}(X;\lag),$ $\WF'(A) =\emptyset$ if and only if $A
\in \Psitwo^{-\infty,l}(X;\lag).$

\item
If $A\in \Psitwo^{k,l}(X;\lag),$ $B\in \Psitwo^{k',l'}(X;\lag)$ then
$$
\sigmatwo_{k+k'-1,l+l'}(i[A,B]) = \{\sigmatwo_{k,l}(A),\sigmatwo_{k',l'}(B)\}
$$
where the Poisson bracket on the right hand side
is computed with respect to the symplectic form
on $\pssp$ lifted from the symplectic form on $T^*X.$

\item\label{ellipticparametrix}
There is a microlocal parametrix near elliptic points: if $p\in\liptic(A)$,
$A\in\Psitwo^{m,l}(X;\lag)$ then there exist $B\in\Psitwo^{-m,-l}(X;\lag)$,
$E,F\in\Psitwo^{0,0}(X;\lag)$ such
that $p\notin\WF'(E)$, $p\notin\WF'(F)$, and $AB=I+E$, $BA=I+F$,
where\footnote{Note that the (non-)vanishing of $(\sbdf^{l-m}\sigmatwo_{m,l}(A))(p)$
is independent of the choice of the defining function $\sbdf$ of the front
face of $\pssp$, so one may reasonably write $\sigmatwo_{m,l}(A)(p)\neq 0$,
meaning $(\sbdf^{l-m}\sigmatwo_{m,l}(A))(p)\neq 0$.}
$$
\liptic (A) = \{p\in\pssp:\ (\sbdf^{l-m}\sigmatwo_{m,l}(A))(p)\neq 0\}
\subset \pssp.
$$

\item\label{mappingproperties}
If $A\in\Psitwo^{m,m}(X;\lag)$ then $A:h^kL^2\to h^{k-m}L^2$ for all $k$.

If $A \in \Psitwo^{-\infty,l}(X, \lag)$ then
$$
A: L^2(X) \to I_{(-l)}^\infty(\lag),
$$
where, for $k \in \NN,$
$$
I_{(s)}^k(\lag) = \{u: h^{-j-s} A_1\dots A_j u \in L^2\  \forall A_i \in
\Psih^1(X),\ \sigma(A_i)\restrictedto_\lag =0, j\leq k\}
$$ (hence $I_{(-l)}^\infty(\lag)$ is, by definition, a space of
semiclassical Lagrangian distributions).  For general $k \in \RR,$
$I_{(s)}^k(\lag)$ is defined by interpolation and duality.

More generally,
if $A \in \Psitwo^{m,l}(X, \lag)$ and $m,k \in \NN,$ then
$$
A: I^k_{(s)}(\lag) \to I^{k-m}_{(s-l)}(\lag).
$$
for each $k.$

The distributions in $I^{-\infty}_{(s)}=\bigcup_k I^{k}_{(s)}$, are called
non-focusing relative to $h^{-s}L^2$ at $\lag$ (of order $-k$, if they are
in $I^{k}_{(s)}$) in \cite{MVW1}.

\item
For $u\in I^{-\infty}_{(l)}$,
there is an associated wavefront set, $\WFtwo^{m,l} u \subset\pssp,$
defined by
$$
(\WF^{m,l} u)^c = \bigcup \{ \liptic(A): A \in \Psitwo^{m,l}, Au \in L^2\}.
$$

$\WF^{\infty,l}(u)=\emptyset$ if and only if $h^{-l} u$ is an
$L^2$-based semiclassical Lagrangian distribution with respect to $\lag$
(as defined above in \ref{mappingproperties}).

Moreover, if $A\in\Psitwo^{m',l'}(X;\lag)$ then $$\WF^{m-m',l-l'}(Au)
\subset\WF^{m,l}(u).$$

Away from $SN(\lag),$ this wavefront set just reduces to the usual
semiclassical wavefront set:
$$
\WFtwo^{m,l} (u) \backslash SN(\lag) = \WF_h^m (u)\backslash \lag,
$$
where we have identified the complement of the front face of $[T^*X;\lag]$
with $T^*X \backslash \lag$ in the natural way.

\item\label{hcalcintwocalc}
The smoothing semiclassical calculus lies inside $\Psitwo(X;\lag):$
$\Psism^k(X) \subset \Psitwo^{k,k}(X;\lag);$ if $A \in \Psism^k(X),$
$$
\sigmatwo_{k,k}(A) = \beta^*\sigma_h(A),\ \WF'(A)=\beta^{-1}(\WF'_h(A)),
$$ (hence $\sigmatwo_{k,k}(A)$
is independent of the fiber variables of $SN(\lag)$); here
$\beta: [T^*X; \lag] \to T^*X$ is the blowdown map, $\sigma_h$ is the
semiclassical principal symbol, and $\WF'_h$ is the semiclassical
operator wave front set.

If $A \in \Psism^k(X)$ and $\sigma_h(A)=0$ on $\lag,$ then
$$
A \in \Psitwo^{k,k-1}(X).
$$

\item\label{awayfromlag}
If $Q\in\Psism^{m'}(X)$ and $\WF'_h(Q)\cap \lag=\emptyset$ then
for all $A\in\Psitwo^{m,l}(X;\lag)$, $QA,AQ\in\Psism^{m+m'}(X)$,
$\WF'(QA),\WF'(AQ)\subset\WF'(Q)$, i.e.\ microlocally away from $\lag$,
$\Psitwo(X;\lag)$ is just $\Psism(X)$.

\item
If $Q,Q'\in\Psism^{0}(X)$ and $\WF'(Q)\cap\WF'(Q')=\emptyset$,
then for $A\in\Psitwo^{m,l}(X;\lag)$, $QAQ'\in\Psitwo^{-\infty,l}(X;\lag)$,
i.e.\ $\Psitwo
(X;\lag)$ is 2-microlocal, but not microlocal at $\lag$: \emph{Lagrangian}
singularities can spread along $\lag$.

\end{enumerate}
\renewcommand{\theenumi}{\Roman{enumi}}

The most important case is the model case, where $\lag$ is the zero section
of $T^*X$. We give the detailed construction arguments in this case: the
definition is in Definition~\ref{def:two-micro-calc}, while the precise
location of the proofs of the properties listed above is given after the
proof of Lemma~\ref{lemma:off-lag}. The general definition is given in
Definition~\ref{def:two-micro-lag}, and the properties are briefly
discussed afterwards.

\begin{remark}\label{remark:noncompact}
While
we chose to exclude diagonal singularities for elements of
$\Psitwo(X;\lag)$ because this is irrelevant for most considerations here,
and because it would require an additional filtration, principal symbol, etc.,
the properties listed easily allow one to define a new space of operators,
\begin{equation}\label{eq:big-calc}
\Psitwo^{m,k,l}(X;\lag)=\Psih^{m,k}(X)+\Psitwo^{k,l}(X;\lag),
\end{equation}
and deduce the analogues of all listed properties. In particular,
note that if $A\in\Psih^{m,k'}(X)$ and $B\in\Psitwo^{k,l}(X;\lag)$ then
choosing $Q\in\Psism^0(X)\subset\Psih^{-\infty,0}(X)$ with $\WF'(\Id-Q)
\cap\lag=\emptyset$ and $\WF'(\Id-Q)\cap\WF'(B)=\emptyset$
then $B=QB+(\Id-Q)B$, $(\Id-Q)B=B-QB\in\Psih^{-\infty,-\infty}(X)$ in view
of the composition formula,
so $A(\Id-Q)B\in\Psih^{-\infty,-\infty}(X)$, while $AQB=(AQ)B
\in\Psitwo^{k+k',l+k'}(X)$ since $AQ\in\Psism^{-\infty,k'}(X)$, so we deduce
that $AB\in\Psitwo^{k+k',l+k'}(X)$.
\end{remark}

\section{The Model Case and the Construction of the Calculus}

We construct $\Psitwo(X;\lag)$ by constructing it first in the model case
of the zero section in $\RR^n$ (i.e.\ $\lag = o\subset T^*\RR^n$) and verifying its
properties, concluding with invariance under semiclassical FIOs preserving
the zero section.


Recall that in the case at hand, our ``total symbol space'' is defined as
$$
\symbsp = [(T^*\RR^n)\times [0,1); o\times 0],
$$
while the ``principal symbol space'' is the side face (the lift of
$T^*\RR^n\times 0$), which can be identified with
$$
\pssp=[T^*\RR^n;o].
$$
Let $\bdf_\sidef$ and $\bdf_\ff$ denote boundary defining functions for
the side and front faces of this blown-up space.
The space of symbols with which we will be primarily concerned will be
$$
S^{m,l}(\symbsp) = \bdf_\sidef^{-m} \bdf_\ff^{-l} \CI_c(\symbsp)
$$
It is also sometimes useful to consider
the space of symbols which are Schwartz
at `fiber infinity',
$$
\dot S^{m,l}(\symbsp) = \bdf_\sidef^{-m} \bdf_\ff^{-l} \cS(\symbsp);
$$
here $\cS(\symbsp)$ stands for the space of Schwartz functions on $\symbsp$,
i.e.\ elements of $\CI(\symbsp)$, which near infinity in $T^*X\times[0,1)$
(where the blow-up of the zero section can be ignored) decay rapidly
together with all derivatives corresponding to the vector bundle structure
(recall that Schwartz functions on a vector bundle are well-defined).

Explicitly, in local coordinates $x=(x_1,\ldots,x_n)$ in some open set $\U
\subset X$ and canonical dual coordinates $\xi=(\xi_1,\ldots,\xi_n)$,
coordinates on $T^*X\times [0,1)_h$ are given by $x,\xi,h$, and $o\times 0$ is
given by $\xi=0$, $h=0.$  Coordinates on $\symbsp$ near the corner (given
by the intersection of the front face with the lift of the boundary,
$h=0$), where $|\xi_k|>\ep|\xi_j|$ for $j\neq k$, are given by
$x,h/|\xi_k|,|\xi_k|$ and $\xi_j/|\xi_k|$ ($j\neq k$), while $x,\Xi=\xi/h$
and $h$ are valid coordinates in a neighborhood of the interior of the
front face. Alternatively, near the corner, one can use polar coordinates,
$x,h/|\xi|,|\xi|$ and $\xi/|\xi| \in \sphere^{n-1}$. Locally $|\xi|$ is
then a defining function for $\ff$, and $h/|\xi|$ is a defining function
for $\sidef$.  Thus, a typical example of an element of $S^{m,l}(\symbsp)$
is a function of the form $h^{-m}|\xi|^{-l+m}a(x,h/|\xi|,|\xi|,\xi/|\xi|)$,
$a\in\CI_c(\U\times [0,\infty)\times [0,\infty)\times\sphere^{n-1})$.

Slightly more globally in the fibers of the cotangent bundle
(but locally in $\U$), one can
use $\langle \xi/h\rangle^{-1}=(h/\xi)(1+(h/|\xi|)^2)^{-1/2}$
as the defining function for $\sidef$, which is now a non-vanishing smooth
function in the interior of $\ff$, so $h\langle \xi/h\rangle$ can be
taken as the defining function of $\ff$. A straightforward calculation
shows that
$a\in \dot S^{-\infty,l}(\symbsp)$ if and only if $b(x,\Xi,h)=h^{l}a(x,h\Xi,h)
\in\CI(\U\times [0,1)_h;\cS(\RR^n_{\Xi}))$, with $\cS$ standing for the
space of Schwartz functions. Indeed, this merely requires noting
that $\Xi_j\pa_{\Xi_k}=\xi_j\pa_{\xi_k}$, and the rapid decay in $\Xi$
fibers corresponds bounds by $C_N\langle \Xi\rangle^{-N}
=C_N\langle\xi/h\rangle^{-N}=C_N\bdf_\sidef^N$.

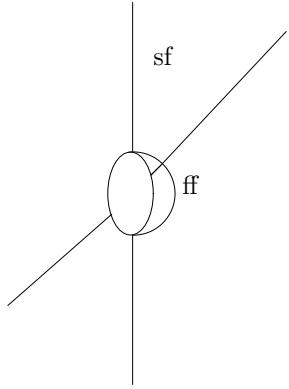
\begin{figure}
\input{symbolspace.customized.eepic}
\caption{The total symbol space $\symbsp,$ in the case $n=2$ and $\lag=0,$
  with base variables omitted.  The front face of the blowup is labeled ff.
  The side face, labeled sf, is the space $\pssp$ on which principal
  symbols are defined, and is canonically diffeomorphic to $[T^*X; o]$  The
  boundary sphere of this side face is diffeomorphic to $SN(o).$}\end{figure}

Let $\Opl,\Opr$ and $\Opw$ denote left-, right-, and Weyl-semiclassical
quantization maps on $\RR^n$, i.e.\ for $a\in S^{m,l}(\symbsp),$
\begin{equation*}\begin{split}
\Opl (a) &= \frac{1}{(2\pi h)^n} \int e^{i (x-y)\cdot\xi/h}
\chi(x,y) a(x,\xi,h) \, d \xi,\\
\Opr (a) &= \frac{1}{(2\pi h)^n} \int e^{i (x-y)\cdot\xi/h}
\chi(x,y) a(y,\xi,h) \, d \xi,\\
\Opw (a) &= \frac{1}{(2\pi h)^n} \int e^{i (x-y)\cdot\xi/h}
\chi(x,y) a((x+y)/2,\xi,h) \, d \xi,
\end{split}\end{equation*}
where $\chi$ is a cutoff properly supported near the diagonal (used to
obtain proper supports), identically $1$ in a smaller neighborhood of the
diagonal, e.g.\ $\chi=\chi_0(|x-y|^2)$, $\chi_0\in\CI_c(\RR)$ identically
$1$ near $0$.
Note that the allowed singularity of $a$ at
$\xi=h=0$ does not cause any problem in defining the integral for $h>0.$
More generally, if
$$
a\in S^{m,l}([\RR^{2n}_{xy}\times\RR^n_\xi\times[0,1)_h;
\RR^{2n}\times\{0\}\times\{0\}])=\CI(\RR^n_x;S^{m,l}(\symbsp_{y,\xi,h})),
$$
we write
$$
I(a)=\frac{1}{(2\pi h)^n} \int e^{i (x-y)\cdot\xi/h}
\chi(x,y) a(x,y,\xi,h) \, d \xi.
$$

\begin{definition}
If $a \in S^{m,l}(\symbsp),$ let
$\esssupp(a)=\esssupp_l(a)$ be the subset of $\pssp$ defined as follows:
$$
\esssupp(a)^c=\{p\in\pssp:\ \exists \phi\in\CI(\symbsp),\ \phi(p)\neq 0,
\ \phi a\in\bigcap_{m'\in\RR} S^{m',l}(\symbsp)\},
$$
i.e.\ a point $p$ is {\em not} in $\esssupp(a)$ is $p$ has a neighborhood
{\em in} $\symbsp$ in which $a$ vanishes to infinite order {\em at}
$\pssp$. We usually suppress the subscript $l$ in the notation.
\end{definition}

We give a manifestly invariant definition of the residual operators in our
calculus: they are powers of $h$ times families of smoothing operators,
with conormal regularity in $h:$

\begin{definition}\label{def:residual}
For each $l \in \RR,$ let
$$\residual^l= \{R: \CmI(\RR^n) \to \CI(\RR^n):
\norm{h^l (h\pa_h)^\alpha\ang{\Lap}^\beta R\ang{\Lap}^\gamma u}\ \leq C_{\alpha\beta\gamma}\norm{u} \
\forall \alpha,\beta,\gamma\}.$$  Let $\residual = \bigcup_{l\in \RR}
\residual^l.$  We further assume that all operators in $\residual$ have
properly supported Schwartz kernels.  (Norms are with respect to $L^2$.)  
\end{definition}

An alternate characterization is as follows.  We let $\kappa(\cdot)$ denote
the Schwartz kernel of an operator.
\begin{lemma}\label{lemma:residual}
$R \in \residual^l$ if and only if
\begin{equation}\label{schkernel}
\sup \abs{h^l (h\pa_h)^\alpha \pa_x^\beta \pa_y^\gamma \kappa(R) (x,y,h)}\leq C_{\alpha\beta\gamma}.
\end{equation}
\end{lemma}
\begin{proof}
Certainly if \eqref{schkernel} does hold, we obtain a uniform estimate on
operator norms as required by Definition~\ref{def:residual}, as indeed we
may estimate Hilbert-Schmidt norms of $R$ in terms of the estimates
\eqref{schkernel} and the size of the support.  Conversely,
Definition~\ref{def:residual} tells us that $h^l (h\pa_h)^\alpha R: H^{-s} \to
H^{s+\abs{\beta}}$ for any desired $s\in \RR$ and multiindex $\beta;$
taking $s>n/2$ and using Sobolev imbedding gives $\pa_x^\beta\pa_y^\gamma
h^l (h\pa_h)^\alpha \kappa(R) \in L^\infty$ (uniformly in $h$).
\end{proof}

As mentioned above
this definition generalizes immediately to a manifold $X$ without boundary:
\begin{definition}\label{def:residual-mfld}
$R\in\residual^l(X)$ if it has a properly supported Schwartz kernel
on $X\times X\times [0,1)$, satisfying \eqref{schkernel} in local coordinates
$x$, resp.\ $y$, or equivalently, that for all $k$ and all compactly
supported vector fields
$V_1,\ldots,V_k\in\Vf(X\times X\times[0,1))$ tangent to $h=0$, there exists
$C$ such that
$|h^l V_1\ldots V_k \kappa(R)|\leq C$.
\end{definition}

Returning to $\RR^n$,
we now show that the quantizations of the ``residual'' symbols
in 
$$
S^{-\infty,l}(\symbsp) \equiv \bigcap_{m\in \RR} S^{m,l}(\symbsp)
$$
lie in the space of residual operators $\residual^l.$

\begin{lemma}\label{lemma:q-smoothing}
$\Opl,\Opr,\Opw$ map $S^{-\infty}(\symbsp)$ into $\residual,$ and
$I$ maps $\CI(\RR^n;S^{-\infty}(\symbsp))$ into $\residual.$
\end{lemma}

\begin{proof}
We have
$$
S^{-\infty,l}(\symbsp) = \{a: a \in \bdf_\sidef^{\infty}\bdf_\ff^{-l}
\CI_c(\symbsp)\}.
$$ Then the kernel of the quantization of $a$ is given by
$$
\Opl(a) = \frac{1}{(2\pi h)^n} \int e^{i(x-y)\cdot\xi/h}
\chi(x,y) a(x,\xi,h) \, d \xi.
$$ Owing to its rapid vanishing at $\sidef,$ we find that $a$ is classical
conormal (i.e.\ a power of a boundary defining function times a $\CI$
function) on the space obtained from $\symbsp$ by blowing down $\sidef,$
i.e.\ by introducing new variables $\Xi = \xi/h$ instead of $\xi;$
we can write $a(x,h\Xi,h)=\tilde a(x,\Xi,h) h^{-l}$ where $\tilde a$ is
$\CI$ and vanishing rapidly as $\Xi \to \infty.$ Hence
\begin{align*}
\Opl(a) &= 
\frac{1}{(2\pi h)^n} \int e^{i(x-y)\cdot\xi/h}
\chi(x,y) a(x,\xi,h) \, d \xi\\ &= 
\frac{h^{-l}}{(2\pi)^n} \int e^{i(x-y)\cdot\Xi}
\chi(x,y) \tilde{a}(x,\Xi,h) \, d \Xi=h^{-l}\chi(x,y)
(\fcal^{-1}\tilde{a})(x,x-y,h),
\end{align*}
where $\fcal^{-1}$ is the inverse Fourier transform in the second argument
of $\tilde a$ (i.e.\ in $\Xi$),
and this is just $h^{-l}$ times a family of smoothing operators with
parameter $h.$ In particular \eqref{schkernel} is easily verified, hence
$\Opl(a)$ lies in $\residual^l.$ Analogous arguments hold for $\Opr,\Opw$
and $I$.
\end{proof}

In fact, we have the following slight strengthening:

\begin{lemma}\label{lemma:q-smoothing-imp}
If $A\in\Op_\bullet(S^{m,l}(\symbsp))$,
and for all $N$ there exists $a_N\in S^{m-N,l}$ such that
$A=\Op_\bullet(a_N)$, then $A\in\residual^l$, where $\bullet$ can be l,r,
or W.
\end{lemma}

\begin{proof}
We prove the lemma for the case of $\Opl;$ the other cases are analogous.

As
\begin{equation*}\begin{split}
&(h\pa_h)^\alpha \pa_x^\beta\pa_y^\gamma\Opl(a_N)\\
&= 
\frac{1}{(2\pi h)^n} \int e^{i(x-y)\cdot\xi/h}
(h\pa_h-i(x-y)\cdot\xi/h)^{\alpha}
(\pa_x+i\xi/h)^\beta\\
&\qquad\qquad\qquad(\pa_y-i\xi/h)^\gamma(\chi(x,y) a_N(x,\xi,h)) \, d \xi,
\end{split}\end{equation*}
the conclusion follows from choosing $N$ large enough such that
$\langle \xi/h\rangle^{|\alpha|+|\beta|+|\gamma|}a_N\in S^{l,l}(\symbsp)$,
which in turn is possible as $\langle \xi/h\rangle=(1+|\xi/h|^2)^{1/2}$
is the reciprocal of a defining function of $\sidef$ as described at the
beginning of the section.
\end{proof}

\begin{lemma}\label{lemma:leftright}
If $a\in\CI(\RR^n;S^{m,l} (\symbsp))$, then there
exist $a_l,a_r,a_W\in S^{m,l}(\symbsp)$ and $R_l,R_r,R_W\in\residual^l$
such that
$$
I(a)=\Opl(a_l)+R_l=\Opr(a_r)+R_r=\Opw(a_W)+R_W,
$$
and
$$
a|_{y=x}-a_l, a|_{x=y}-a_r,a|_{x=y}-a_W\in S^{m-1,l}(\symbsp).
$$ In particular, we may change from left- to right-quantization and
vice-versa: if $a \in S^{m,l} (\symbsp);$ then there exist $b,b' \in
S^{m,l}(\symbsp)$ and $R,R' \in \residual^l$ such that
$$
\Opl(a) = \Opr(b)+R,\quad \Opr(a) = \Opl(b')+R'.
$$

Moreover, $\esssupp a_l=\esssupp a_r=\esssupp a_W$.
\end{lemma}

\begin{proof}
To begin, we prove that for $a\in\CI(\RR^n;S^{m,l} (\symbsp))$, $I(a)=\Opr(b)+R$,
$b\in S^{m,l}(\symbsp)$, $R\in\residual^l$. The statement about
$\Opl$ follows the same way reversing the role of $x$ and $y$ below.

Using a partition of unity, we may decompose $a$ into pieces supported on
the lifts to $\symbsp$ of the set $\{\xi_j\neq 0\}\subset (T^*\RR^n \times [0,1))$ for
various values of $j.$ By symmetry, it will suffice to deal with the term
supported on $\xi_1 \neq 0.$ On this region, we may take as coordinates in
$\symbsp$ the functions $\Xi = \xi'/\xi_1,$ $H= h/\xi_1$, and $\xi_1.$
Thus, $\xi_1$ is locally a defining function for $\ff$ and $H$ for
$\sidef.$ We may Taylor expand in $x$ around $y$,
$$
a(x,y,\xi,h) \sim \sum \frac{1}{\alpha!} (x-y)^\alpha
(\pa_x^\alpha a)(y,y,\xi,h))
$$
Now in the variables $\Xi,\xi_1,H,$ we
have
\begin{equation}\label{vflifts}
\begin{aligned}
h \pa_{\xi_1} &= H(\xi_1 \pa_{\xi_1}-H\pa_H-\Xi\cdot \pa_\Xi),\\
h \pa_{\xi'} &= H \pa_\Xi,\\
h \pa_h &= H\pa_H,
\end{aligned}
\end{equation}
hence by our symbolic assumptions on $a,$
$$
(h\pa_\xi)^\alpha (\pa_x^\alpha a)(x,x,\xi,h))\in H^{-m+\abs{\alpha}}
\xi_1^{-l} \CI(\RR^n \times\symbsp)
$$ near the corner $H=\xi_1=0,$ hence these terms may be Borel summed to
some $a_r\in H^{-m} \xi_1^{-l}\CI(\RR^n \times \symbsp).$

For any $N\in \NN,$ by integrating by parts, we have
\begin{align}
I(a) &=\sum_{\abs{\alpha}<N} \frac{1}{(2\pi h)^n} \int e^{i (x-y)\cdot \xi/h}
\chi(x,y) \frac{1}{\alpha!}(x-y)^\alpha  (\pa_x^\alpha a)(y,y,\xi,h)\, d\xi+R_N'\\
&= 
\sum_{\abs{\alpha}<N} C_\alpha \frac{1}{(2\pi h)^n} \int e^{i (x-y)\cdot \xi/h}
\chi(x,y) \frac{1}{\alpha!}
((h\pa_\xi)^\alpha \pa_x^\alpha a)(y,y,\xi,h)\, d\xi+R_N'
\\&= \Opr(a_r)+R_N
\end{align}
where $R_N$ and $R_N'$ both have the form
\begin{multline}
I\left(\sum_{\abs{\alpha}=N} (x-y)^\alpha r'_\alpha\right)\\
= 
\sum_{\abs{\alpha}=N} \frac{1}{(2\pi h)^n} \int e^{i (x-y)\cdot \xi/h}
\chi(x,y) (h\pa_\xi)^\alpha (r'_\alpha(x,y,\xi,h))\, d\xi
\end{multline}
with $r'_\alpha\in \CI(\RR^n \times S^{-N,l}(\symbsp)).$  We thus have
$$
I(a)-\Opr(a_r) \in I(\CI(\RR^n;S^{-N,l}(\symbsp)))
$$ for all $N \in \NN.$ By
Lemma~\ref{lemma:q-smoothing-imp} we obtain the desired result.

The statement about $\Opw$ can be proved similarly, writing
$$
a(x,y,\xi,h)=\tilde a((x+y)/2,(x-y)/2,\xi,h),
$$
i.e.\ $\tilde a(w,z,\xi,h)
=a(z+w,z-w,\xi,h)$, and expanding $\tilde a$ in Taylor series in $z=(x-y)/2$
around $0$, so
$$
a(x,y,\xi,h)\sim \sum \frac{1}{\alpha!}\left(\frac{x-y}{2}\right)^\alpha
\left((\pa_x+\pa_y)^\alpha a\right)\left(\frac{x+y}{2},\frac{x+y}{2},\xi,h\right).
$$

Finally, the statements about $\esssupp a_l$, etc., follows for e.g.\ if
$a(x,y,\xi,h)=a_l(y,\xi,h)$, the terms
$a_{r,\alpha}=\frac{1}{\alpha!}
(h\pa_\xi)^\alpha \pa_x^\alpha(a_l)(y,\xi,h)$ in the asymptotic
expansion for $a_r$ all have $\esssupp a_{r,\alpha}\subset \esssupp a_l$.
\end{proof}

\begin{lemma}\label{lemma:off-diag}
If $a\in \CI(\RR^n;S^{m,l}(\symbsp))$
then for $\phi\in\CI(X^2\times[0,1))$ with support
disjoint from $\diag\times\{0\}$, $\phi I(a)\in\residual^l$.
\end{lemma}

\begin{proof}
As $I(a)\in\CI(X^2\times(0,1))$, we may assume that $\supp\phi$
is disjoint from $\diag\times[0,1)$, hence $|x-y|>\ep$ on $\supp\phi$.
Then for all $N$,
\begin{equation*}\begin{split}
I(a)&=\frac{1}{(2\pi h)^n} \int \frac{h^{2N}}{|x-y|^{2N}}
\Delta_\xi^{N}e^{i (x-y)\cdot\xi/h}
\chi(x,y) a(x,y,\xi,h) \, d \xi\\
&=\frac{1}{(2\pi h)^n} \int \frac{1}{|x-y|^{2N}}
e^{i (x-y)\cdot\xi/h}
\chi(x,y) (h^{2N}\Delta_\xi^N a)(x,y,\xi,h) \, d \xi
\end{split}\end{equation*}
We can assume, using a partition of unity as above, that $a$ is supported
in the lift of the set where $\xi_1\neq 0$. Then \eqref{vflifts}
shows that $h^{2N}\Delta_{\xi}^N a
\in S^{m-2N,l}(\symbsp)$.
Choosing $N$ sufficiently large, depending
on $\alpha,\beta,\gamma$, it follows
immediately (cf.\ the proof of Lemma~\ref{lemma:q-smoothing-imp}) that
$$
|h^l (h\pa h)^{\alpha}\pa_x^\beta \pa_y^\gamma I(a)|\leq C,
$$
for each derivative at most gives an additional factor of $H^{-1}$ in growth.
\end{proof}

The proof of this lemma can in fact be extended to show that $\Opl(a)$
determines $a$ modulo $S^{-\infty,l}(\symbsp):$

\begin{lemma}\label{lemma:Opl-to-a}
For $a\in S^{m,l}(\symbsp)$, let $\kappa$ denote the Schwartz kernel
of $\Opl(a)$. Then
$$
a(x,\xi,h)-(\fcal_z\kappa(x,x-z,h))(\xi)\in \dot S^{-\infty, l}(\symbsp),
$$
where
$\dot S$ denotes the space of symbols rapidly decreasing at infinity
(rather than compactly supported)
in $\symbsp$.
In particular, modulo $S^{-\infty,l}$, $\Opl(a)$ determines $a$.\footnote{Here
$\tilde a(x,\xi,h)=(\fcal_z\kappa(x,x-z,h))(\xi)\in\CI(T^*\RR^n\times(0,1))$,
so one may indeed think
of it as a function, as indicated by the notation. However, one must
use the blown-up coordinates on $[T^*\RR^n\times[0,1);o\times 0]$ in order
to realize $a$ as a polyhomogeneous function, hence in order to evaluate
$\bdf_{\sidef}^m\bdf_{\ff}^l\tilde a$ at the front face.}

Analogous statements hold for $\Opr(a)$ and $\Opw(a)$ as well.
\end{lemma}

\begin{proof}
For $a\in S^{m,l}(\symbsp)$, let
$$
K(x,z,h)=(\fcalh^{-1}_{\xi}a(x,\xi,h))(z)\equiv (2\pi h)^{-n}\int
e^{i z\cdot\xi/h}a(x,\xi,h)\,d\xi
$$
be the semiclassical inverse Fourier transform of $a$ in $\xi$,
so
$$
a(x,\xi,h)=(\fcalh_{z} K(x,z,h))(\xi)=(\fcal_z K(x,z,h))(\xi/h)
=\int e^{-iz\cdot\xi/h} K(x,z,h)\,dz.
$$
Then
$$
\kappa(x,x-z,h)=\chi(x,x-z) K(x,z,h),
$$
hence
$$
r \equiv a(x,\xi,h)-(\fcal_z\kappa(x,x-z,h))(\xi) = \fcal_z\left(
(1-\chi(x,x-z))K(x,z,h)\right)(\xi/h);
$$
we need to show that this lies in $S^{-\infty,l}(\symbsp).$

The proof of the preceding lemma shows that $(1-\chi(x,y))K(x,y,h)$ is
Schwartz in $x-y$, smooth in $x$, conormal in $h$ of order $l$, i.e.
$$
h^{l}(h\pa_h)^s z^\alpha D_z^\beta D_x^\gamma
(1-\chi(x,x-z)) K(x,z,h)
$$
is bounded for all $s,\alpha,\beta,\gamma$, so its (non-semiclassical)
Fourier transform
in $z=x-y$,
\begin{equation*}\begin{split}
\tilde r(x,\Xi,h)&=\left(\fcal_z (1-\chi(x,x-z)) K(x,z,h)\right)(\Xi)\\
&=\int e^{-i\Xi\cdot z}  (1-\chi(x,x-z)) K(x,z,h)\,dz
\end{split}\end{equation*}
satisfies
$$
h^{l}(h\pa_h)^s \Xi^\alpha D_\Xi^\beta D_x^\gamma
\tilde r(x,\Xi,h)\in L^\infty
$$
for all $s,\alpha,\beta,\gamma$. Thus, $\tilde r\in \CI(\RR^n_x\times[0,1)_h;
\cS(\RR^n_\Xi))$, so as remarked at the beginning of the section,
$$
r(x,\xi,h)=\tilde r(x,\xi/h,h)\in \dot S^{-\infty,l}(\symbsp).\qed
$$
\noqed
\end{proof}

We now prove diffeomorphism invariance.
\begin{lemma}\label{lemma:diffeo-inv}
If $F:U\to U'$ is a diffeomorphism, $G=F^{-1}$,
$U,U'\subset\RR^n$,
and $A=\Opl(a)$, $a\in S^{m,l}(\symbsp)$, with the Schwartz kernel of
$A$ supported in $U'$ then $F^*AG^*=\Opl(b)+R$ for some
$b\in S^{m,l}(\symbsp)$, $R\in \residual^l$.

Moreover, $b-(G^\sharp)^*a\in S^{m-1,l}(\symbsp)$, where $G^\sharp:T^*U
\to T^*U'$ is the induced pull-back of one-forms, and
$\esssupp b=G^\sharp(\esssupp a)$.

For the Weyl quantization, and with $A$ acting on half-densities,
the analogous statement holds with the
improvement $b-(G^\sharp)^*a\in S^{m-2,l}(\symbsp)$.
\end{lemma}

\begin{proof}
We follow the usual proof of the diffeomorphism invariance formula.

Note first that $\residual^l$ is certainly invariant under pullbacks by
diffeomorphisms, and a partition of unity, with an element identically
$1$ near the diagonal, allows us to assume that $K_A$ is supported in
a prescribed neighborhood of the diagonal.

The Schwartz kernel $K_B(x,y) |dy|$ of $B=F^*AG^*$ is $(F\times F)^*
(K_A(\tilde x,\tilde y) |d\tilde y|)$, i.e.
\begin{equation*}\begin{split}
K_B(x,y)
&=K_A(F(x),F(y))|\det F(y)|\\
&=(2\pi h)^{-n}\int e^{i(F(x)-F(y))\cdot \xi/h}
a(F(x),\xi)|\det F(y)|\,d\xi.
\end{split}\end{equation*}
Now, $F_i(x)-F_i(y)=\sum_j (x_j-y_j)T_{ij}(x,y)=T(x,y)(x-y)$
by Taylor's theorem,
with $T_{ij}(x,x)=\pa_j F_i(x)$, so $T(x,x)$ invertible. Thus,
$T$ is invertible in a neighborhood of the diagonal; we take this as
the prescribed neighborhood mentioned above. Then, with $\eta=T^t(x,y)\xi$,
\begin{equation*}\begin{split}
&K_B(x,y)\\
&=(2\pi h)^{-n}\int e^{i(x-y)\cdot \eta/h}
a(F(x),(T^t)^{-1}(x,y)\eta)|\det T(x,y)|^{-1}|\det F(y)|\,d\eta.
\end{split}\end{equation*}
By Lemma~\ref{lemma:leftright}, this is of the form $\Opl(b)+R'$ for some
$b\in S^{m,l}(\symbsp)$, $R'\in \residual^l$, provided that we show that
$a(\dots)\in\CI(\RR^n;S^{m,l}(\symbsp))$, which in turn is immediate.

For the Weyl quantization, acting on half-densities, the Schwartz kernel
$$
K_B(x,y) |dx|^{1/2}|dy|^{1/2}
$$
of $B=F^*AG^*$ is $(F\times F)^*
(K_A(\tilde x,\tilde y) |d\tilde x|^{1/2}|d\tilde y|^{1/2})$,
$$
K_B(x,y)
=K_A(F(x),F(y))|\det F(x)|^{1/2}|\det F(y)|^{1/2}.
$$
Now we
use Taylor's theorem for $F$ around $(x+y)/2$,
so $F_i(x)-F_i(y)=\sum_{ij} (x_j-y_j)T_{ij}(x,y)$ with
$T_{ij}(x,y)-T_{ij}((x+y)/2,(x+y)/2)=O(|x-y|^2)$,
and $(F(x)+F(y))/2=F((x+y)/2)+O(|x-y|^2)$, with an analogous statement
for the product of the determinants, to obtain the improved result.
\end{proof}

In view of the diffeomorphism invariance and Lemma~\ref{lemma:off-diag},
we can naturally define 2-microlocal
operators on manifolds, associated to the 0-section.

\begin{definition}\label{def:two-micro-calc}
Let $$\Psitwo^{m,l}(\RR^n;o)= \{\Opl(a): a \in S^{m,l}(\symbsp)\} + \residual^l.$$
If $X$ is a manifold without boundary, let $\Psitwo^{m,l}(X,o)$
consist of operators $A$ with properly supported Schwartz kernels\footnote
{Here $\pi_R:X\times
X\times[0,1)\to X$ is the
projection to the second factor of $X$, $\Omega X$ the
density bundle.}
$K_A\in\CmI(X\times X\times[0,1);\pi_R^*\Omega X)$,
such that for any coordinate neighborhood $U$ of $p\in X$, and any
$\phi,\psi\in\CI_c(
U),$ $A$ satisfies $\phi A\psi\in\Psitwo^{m,l}(\RR^n;o)$, while if
$\phi,\psi\in\CI_c(X)$ with disjoint support then
$\phi A\psi\in\residual^l(X)$.
\end{definition}

\begin{remark}\label{remark:semicl}
Directly from the definition,
$$
\Psism^m(\RR^n)=h^{-m}\Psism^0(\RR^n)
\subset\Psitwo^{m,m}(\RR^n),
$$
with the relationship between total symbols,
modulo $S^{-\infty,m}(\symbsp)$,
for, say, left-quantization, given by the pullback under the blow-down map
$$
\beta:[T^*\RR^n\times [0,1);o\times 0].
$$
\end{remark}

\begin{lemma}\label{lemma:ideal}
$\residual^l\big(\Psitwo^{m',l'}(\RR^n;o)\big)\subset\residual^{l+l'}$,
$\big(\Psitwo^{m',l'}(\RR^n;o)\big)\residual^l\subset\residual^{l+l'}$.
\end{lemma}

\begin{proof}
It suffices by Lemma~\ref{lemma:leftright} to show that if $a \in
S^{m,l}(\symbsp)$ and $R \in \residual^{l'}$ then 
\begin{equation}\label{toshow}
R \circ \Opr(a) \in \residual^{l+l'},
\end{equation}
as the rest of the statement will follow by taking adjoints.  To show
\eqref{toshow}, we begin by showing that $h^{l+l'} R \circ \Opr(a)$ is
bounded on $L^2$ (uniformly in $h$).  If $m$ is negative, the uniform
$L^2$-boundedness of $h^l\Opr(a)$ follows from
Calder\'on-Vaillancourt,\footnote{Note that we use the fact that $h$ lifts
to $\symbsp$ to be $H\xi_1$ in the local coordinates of the proof of
Lemma~\ref{lemma:leftright}, so that $h^l$ times a symbol in
$S^{m,l}(\symbsp)$ is bounded.} so it suffices to consider the case $m>0.$
In that case, let $k$ be an integer greater than $m.$ We again split $a$ up
into pieces and employ local coordinates as in the proof of
Lemma~\ref{lemma:leftright}; thus, using the fact that
$$
\frac{h}{\xi_1}D_{y_1} e^{i(y-z)\cdot \xi/h} = e^{i(y-z)\cdot \xi/h},
$$
we may write
$$
R \circ \Opr(a) = \frac{1}{(2\pi h)^n}\int K(x,y,h) a(z,\xi,h)
\left( \frac{h}{\xi_1} D_{y_1}\right)^k e^{i(y-z)\cdot \xi/h}\, d\xi\, dy,
$$
where $K$ is the kernel of $R.$  We now integrate by parts in $y_1,$ to
obtain 
$$
R \circ \Opr(a) = \frac{1}{(2\pi h)^n}\int (D_{z_1})^k K(x,y,h) a(z,\xi,h)
\left( \frac{h}{\xi_1} \right)^k e^{i(y-z)\cdot \xi/h}\, d\xi\, dy;
$$ noting that $h^l (h/\xi_1)^k a \in S^{0,0}(\symbsp)$ and that $K$ is
smooth in $z,$ we again obtain $L^2$ boundedness by Calder\'on-Vaillancourt.

To finish the proof, it suffices to show that $h^{l+l'} (h\pa_h)^\alpha
D^\beta R \circ \Opr(a)D^\gamma$ are also $L^2$ bounded.  The follows from stable
regularity of the kernel of $R$ under $h\pa_h$ and $D_x,$ and from a
further integration by parts, since
$$
D_z e^{i(y-z)\cdot \xi/h}=-D_y e^{i(y-z)\cdot \xi/h},
$$
and $y$-derivatives falling on $K$ may also be absorbed without loss.
\end{proof}

\begin{theorem}\label{thm:star-alg}
$\Psitwo(\RR^n;o)$ and $\Psitwo(X,o)$ are bi-filtered $*$-algebras,
with $\residual$ a filtered
two-sided ideal.
\end{theorem}
\begin{proof}
By localization we immediately reduce the general case to $\RR^n$.

That $\Psitwo(\RR^n;o)$ is closed under adjoints follows from our
ability to exchange left and right quantization, as proved above, together
with the fact that the residual calculus is closed under adjoints.

To prove that the calculus is closed under composition, it suffices (using
Lemmas~\ref{lemma:ideal} and \ref{lemma:leftright}) to show that if
we take $a \in S^{m,l}(\symbsp)$ and $b \in S^{m',l'}(\symbsp)$ then $\Opl(a)
\circ \Opr(b) \in \Psitwo^{m+m',l+l'}(\RR^n;o)+\residual^{l+l'}.$  We have
\begin{align*}
\Opl(a) \circ \Opr(b)&= \frac{1}{(2\pi h)^{2n}} \int a(x,\xi,h)b(y,\eta,h)
e^{i (x-w) \cdot \xi/h} e^{i(w-y) \cdot \eta/h}\,d\xi\,d\eta\\
&=
\frac{1}{(2\pi h)^{n}} \int a(x,\xi,h)b(y,\xi,h)
e^{i (x-y) \cdot \xi/h} \,d\xi
\end{align*}
Lemma~\ref{lemma:leftright} permits us to rewrite this expression as a left
quantization of a symbol in $S^{m+m',l+l'}(\symbsp)$ plus a term in
$\residual^{l+l'}.$

The ideal property of $\residual$ is immediate from Lemma~\ref{lemma:ideal}.
\end{proof}

We now discuss the definition and properties of the principal symbol map.
If $A \in \Psitwo^{m,l}(\RR^n;o)$ is given by 
$$
A =\Opl(a) +R,\quad R \in \residual^l,
$$ we define\footnote{As $h$ is a globally well-defined function on
$\symbsp$, we do not need to introduce a line bundle to take care of this
renormalization; this is in contrast with the case when one wishes to define
the usual principal symbol as a function on the cosphere bundle, but a line
bundle appears unavoidably in the definition.}
$$
\sigmatwo(A) = (h^m a)\restrictedto_{\sidef}\in S^{l-m}(\sidef)
=S^{l-m}(\pssp).
$$ As
usual, we may write this in terms of the kernel of $\Opl(a)$ in terms of
Fourier transform:\footnote{The following formula should be interpreted
with a grain of salt: the value of the symbol at $\rho=0$ (where, indeed,
it is of greatest interest) must be obtained from the formula by continuous
extension from the case $\rho>0,$ where it makes sense (and equals the
ordinary semiclassical symbol).}
$$(\bdf_\ff^l a)\restrictedto_{\sidef}(x,\rho,\hat\xi) =
\lim_{h\downarrow 0} \rho^l\FT_{z} (\kappa(\Opl(a))(x,x-z)) (\rho\hat\xi);
$$
here
we have identified $\sidef$ with $\pssp=[T^*\RR^n; o],$ and
$(x,\rho,\hat\xi)$ are coordinates in this space, hence
$\rho=\abs{\xi},\hat\xi=\xi/\abs{\xi};$ we use $\kappa$ to denote the
Schwartz kernel of an operator.

Note that $\sigmatwo(A)$ is not a priori well-defined owing to the
presence of the term in $\residual$ in our definition of the calculus,
but Lemma~\ref{lemma:Opl-to-a} shows that in fact it is. Also, directly
from the definition of $\Psitwo^{m,l}(X,o)$, $\sigmatwo(A)$ can be
defined by localization (i.e.\ considering $\phi A\phi$, $\phi$ identically
$1$ near the point in question) for arbitrary $X$, and is independent of
all choices.

\begin{lemma}\label{lemma:ses}
The principal symbol sequence
\begin{equation*}
0 \to \Psitwo^{k-1,l}(X;o) \to \Psitwo^{k,l}(X;o)
\overset{\sigmatwo_{k,l}}{\to} S^{l-k}(\pssp) \to 0,
\end{equation*}
where the map $\Psitwo^{k-1,l}(X;o) \to \Psitwo^{k,l}(X;o)$ is inclusion,
is exact.
\end{lemma}

\begin{proof}
$\sigmatwo_{k,l}$ is surjective since $\sigmatwo_{k,l}(\Opl(a))=a$.
$\sigmatwo_{k,l}(A)=0$ if $A\in\Psitwo^{k-1,l}(X;o)$ directly from 
the definition. Finally, if $\sigmatwo_{k,l}(A)=0$ for some
$A\in\Psitwo^{k,l}(X;o)$, then $A-\Opl(a)\in\Psitwo^{-\infty,l}(X,o)$
for some $a\in S^{k,l}(\symbsp)$, with $h^{k}a|_{\sidef}=\sigmatwo_{k,l}(A)$,
so if the latter vanishes, then $h^{k}a|_{\sidef}\in\bdf_{\sidef}S^{0,l-k}
(\symbsp)$, hence $a\in S^{k-1,l}(\symbsp)$, giving the conclusion.
\end{proof}

\begin{lemma}\label{lemma:symbhomom}
The principal symbol map is a homomorphism, and
if $A\in \Psitwo^{k,l}(X;o),$ $B\in \Psitwo^{k',l'}(X;o)$ then
$$
\sigmatwo_{k+k'-1,l+l'}(i[A,B]) = \{\sigmatwo_{k,l}(A),\sigmatwo_{k',l'}(B)\}
$$
where the Poisson bracket is computed with respect to the symplectic form
on $[T^*X; o]$ lifted from the symplectic form on $T^*X.$
\end{lemma}
This follows from $[A,B]\in\Psitwo^{k+k'-1,l+l'}(X;o)$ (which in turn
follows from
$$
\sigmatwo_{k+k',l+l'}([A,B])=
\sigmatwo_{k,l}(A)\sigmatwo_{k',l'}(B)-\sigmatwo_{k',l'}(B)\sigmatwo_{k,l}(A)
=0
$$
and the exactness of the principal symbol sequence), with the principal
symbol in $\Psitwo^{k+k'-1,l+l'}(X;o)$ calculated
by continuity from the region $\xi\neq 0$,
where the stated formula follows from a standard argument in the
development of the semiclassical calculus.

We now discuss the definition and properties of the operator wave front set.
If $A \in \Psitwo^{m,l}(\RR^n;o)$ is given by 
$$
A =\Opl(a) +R,\quad a\in S^{m,l},\ R \in \residual^l,
$$ we define
$$
\WF'(A) =\WF'_l(A)= \esssupp_l(a)=\esssupp(a).
$$
Again, this is well-defined by
Lemma~\ref{lemma:Opl-to-a}, and by localization, $\WF'(A)$ is also
well-defined for $A\in\Psitwo^{m,l}(X,o)$; it is a subset of $\pssp.$
Directly from the proof of
Theorem~\ref{thm:star-alg}, which in turn hinges on the asymptotic
expansion given in the proof of Lemma~\ref{lemma:leftright}, we have:

\begin{lemma}\label{lemma:WFp}
For $A\in\Psitwo^{m,l}(\RR^n;o)$, $\WF'(A)=\emptyset$ if and only if
$A\in\Psitwo^{-\infty,l}(\RR^n;o)$.

For $A,B\in\Psitwo(\RR^n;o)$, $\WF'(A+B)\subset\WF'(A)\cup\WF'(B)$ and
$\WF'(AB)\subset\WF'(A)\cap\WF'(B)$.
\end{lemma}

The notion of the operator wave front set also allows us to show that
microlocally away from $o$, $\Psitwo^{m,l}(X,o)$ is the same as
$\Psism^m(X)$. If $Q \in\Psism^m(X),$ we let $\WF'_h(Q)$ denote the usual
semiclassical wave front set, while treating $Q$ as an element of
$\Psitwo(X;o),$ we have $\WF'(Q)=\beta^{-1}(\WF'_h(Q))$ directly from the
definition of $\WF'$ and $\WF'_h$ as essential support.

\begin{lemma}\label{lemma:off-lag}
If $Q\in\Psism^{m'}(X)$ and $\WF'_h(Q)\cap o=\emptyset$ then
for all $A\in\Psitwo^{m,l}(X,o)$, and $QA,\ AQ\in\Psism^{m+m'}(X)$,
$\WF'(QA),\WF'(AQ)\subset\WF'(Q)$.

Thus, the set of operators $Q\in\Psism(X)$ with $\WF'_h(Q)\cap o=\emptyset$
is a two-sided ideal in $\Psitwo(X,o)$.
\end{lemma}

\begin{proof}
It is straightforward to see the conclusion when $A$ is residual,
as
$$(h/\xi_j)D_{x_j} e^{i(x-y)\cdot\xi/h}=e^{i(x-y)\cdot\xi/h},
$$
so
integration by parts in $x$, using that $\xi\neq 0$ on $\esssupp q$ (where $Q
=\Opr(q)$),
shows that, if $K$ is the Schwartz kernel
of $R\in\residual^l$,
$$
R\Opl(q)(z,y)=(2\pi h)^{-n}\int K(z,x)\,e^{i(x-y)\cdot\xi/h} q(y,\xi,h)\,d\xi\,dx
$$
lies in $h^\infty\CI(X^2\times[0,1))$.
Thus, we may assume that $A=\Opl(a)$, and $Q=\Opr(q)$,
$q\in h^{-m'}\CI_c(T^*X\times[0,1))$ with $\esssupp q\cap o=\emptyset$. Then
$AQ=I(b)$, $b(x,y,\xi,h)=a(x,\xi,h)q(y,\xi,h)$, as in the proof
of Theorem~\ref{thm:star-alg}, and
$$
b\in S^{m+m',-\infty}(\symbsp)
\subset S^{m+m'}(T^*X\times[0,1)),
$$
so $I(b)\in\Psism^{m+m'}(X)$.
\end{proof}

We now check that the properties listed in Section~\ref{sec:calculus} hold
for $\Psitwo(X,o)$:

\renewcommand{\theenumi}{\roman{enumi}}
\begin{enumerate}
\item
This is Theorem~\ref{thm:star-alg}, plus the observation that
if $A_j-\Opl(a_j)\in\residual^l$, $a_j\in S^{m-j,l}(\symbsp)$, then one
can Borel sum the $a_j$ to some $a\in S^{m,l}(\symbsp)$ with
$a-\sum_{j=0}^{N-1}a_j\in S^{m-N,l}(\symbsp)$ for all $N$.
\item
This is Lemma~\ref{lemma:symbhomom}, Lemma~\ref{lemma:ses},
and the preceding discussion.
\item
One can take $\Op=\Opl$, for instance.
\item
This is Lemma~\ref{lemma:WFp} and the preceding discussion.
\item
This is Lemma~\ref{lemma:symbhomom}.
\item
This is a standard consequence of (i)-(iv).
\item
The uniform boundedness of $\Psitwo^{m,m}(X,o)$ from $h^kL^2$ to
$h^{k-m}L^2$ follows from the uniform boundedness of $\Psitwo^{0,0}(X,o)$
on $L^2$, which in turn is a consequence of the corresponding property of
$\residual^0$ and of the argument of Calder\'on-Vaillancourt, as noted in
the proof of Lemma~\ref{lemma:ideal}.

For $A\in\Psitwo^{-\infty,l}(X,o)$, $A:L^2(X)\to I^\infty_{(-l)}(o)$
follows from the definition of $\residual^l$ and $h^{-j}A_1\ldots A_j
=(h^{-1}A_1)\ldots (h^{-1}A_j)$, $h^{-1}A_i\in\Psitwo^{1,0}(X,o)$,
so $h^{-j}A_1\ldots A_j A\in\Psitwo^{-\infty,l}(X,o)$, hence maps $L^2(X)$
to $h^{-l}L^2(X)$.

As $\mcal=\{A\in\Psism^1(X):\ \sigma(A)|_o=0\}$ is a (locally) finitely generated
module over $\Psism^0(X)$, with any set of $\CI$ vector fields
spanning $T_pX$ for all $p$ giving a set of generators (for vanishing
of the principal symbol at $o$ means that $A=q_L(h^{-1}a)$, $a|_o=0$,
so $a=\sum a_j\xi_j$ in local coordinates by Taylor's theorem), closed
under commutators, we deduce that for non-negative integers $k$,
$$
I^k_{(0)}(o)=L^\infty((0,1)_h;H^k_{\loc}(X)),
$$ hence (by interpolation and duality) in general the same formula still
holds.  Note also that these spaces are local.
In case $X=\RR^n$, using the `large calculus'
discussed at the end of Section~\ref{sec:calculus}, since
$A=(1+\Delta)^{k/2}\in\Psitwo^{k,k,0}(\RR^n;o)$ is
elliptic, we have produced
an elliptic operator $A\in\Psitwo^{k,k,0}(\RR^n;o)$ such that $u\in
I^k_{(0)}$ implies $Au\in L^2$. The elliptic parametrix construction shows
the converse, so for all $k\geq 0$ (multiplying by $h^s$ if needed)
$$
I^k_{(s)}=\{u\in h^{s}L^2:\ \exists\ \text{elliptic}\ A\in\Psitwo^{k,k,0}(X,o),
\ Au\in h^sL^2\}.
$$ Given $A \in \Psitwo^{k,k,0}(X,o)$ elliptic there exists a parametrix
$B\in\Psitwo^{-k,-k,0}(X,o)$, with $BA-\Id=E\in\Psitwo^{-\infty,-\infty,0}(X,o).$ Thus
if $\tilde A\in\Psitwo^{k-m,k-m,0}(X,o)$ is also elliptic then for any 
$P\in\Psitwo^{m,l}(X,o),$
$$
\tilde APu=\tilde AP(BA-E)u=(\tilde APB)(Au)-(\tilde APE)u,
$$ with $\tilde APB\in\Psitwo^{0,l}(X,o)$, $\tilde
APE\in\Psitwo^{-\infty,l}(X,o)$ both bounded $h^sL^2\to h^{s-l}L^2$.  As a
result, we conclude that $P:I^k_{(s)}(o)\to I^{k-m}_{(s-l)}(o),$ whenever
$k,k-m\geq 0.$ The general case follows by duality arguments.

\item
These are standard consequences of the definition of $\WF^{m,l}$, and
the properties of $\Psitwo(X,o)$, in particular (iv), (vi) and (vii).
\item
See Remark~\ref{remark:semicl}.  The case when $\sigma_h(A)=0$ on $o$ is
evident from the fact that we may write such an operator as $A=\Opl(a)$
with $h^{-k}a$ vanishing on $h=\xi=0,$ hence lifting to be in
$S^{k,k-1}(\symbsp).$
\item
See Lemma~\ref{lemma:off-lag}.
\item
Follows from (iv).
\end{enumerate}
\renewcommand{\theenumi}{\arabic{enumi}}

The following result is necessary in order to transfer our definition from
the model case $\lag=o$ to the case of a general Lagrangian in an invariant
way.

\begin{proposition}\label{prop:FIO-o}
Let $T$ be a properly supported semiclassical FIO with canonical relation
$\Phi$ equal to the identity on the zero section, and with $T$ elliptic on
$\U,$ an open set in $T^*\RR^n;$ let $S$ be a microlocal parametrix for $T$
on $\U.$ Then for $A \in \Psitwo^{m,l}(\RR^n; o)$ with $\WF'(A) \subset
\U,$
$$
T A S \in \Psitwo^{m,l}(\RR^n;o),
$$
with
$$
\sigmatwo (T A S) =(\Phi^{-1})^* \sigmatwo
(A),
$$
and
$$
\WF' (T A S)=\Phi(\WF'(A)).
$$
\end{proposition}
The proof proceeds by deformation to a pseudodifferential operator---cf.\
section~10 of \cite{Evans-Zworski} and \cite{Hitrik-Sjostrand1}.
\begin{proof}
The map $\Phi,$ since it is symplectic, takes $(x,\xi) \to (X, \Xi)$ with
$$
(X_i,\Xi_i) = (x_i+O(\xi),\xi_i+O(\xi^2)),
$$ hence can be parametrized by a \emph{generating function} (see
\cite[\S47A]{Arnold}) $S(x,\Xi) =x\cdot\Xi+\sum \Xi_i \Xi_j \tilde{S}(x,\Xi).$ In
a neighborhood of $o,$ we may thus connect $\Phi$ to the identity map via a
family of symplectomorphisms $\Phi_t$ parametrized by $x\cdot \Xi+t\sum
\Xi_i \Xi_j \tilde{S}(x,\Xi).$ Thus, $\Phi_0=\Id,$ $\Phi_1=\Phi,$ and
$\Phi_t$ fixes the zero section for each $t.$ We connect $T$ to a
semiclassical pseudodifferential operator (microlocally near $o$) via a
family $T_t$ of elliptic semiclassical FIOs given by
\begin{equation}\label{FIO}
T_t = h^{(-n+N)/2}\int e^{i \phi(t,x,y,\theta)/h} b(t,x,y,\theta;h) \,d\theta,
\end{equation} with $T_0\in \Psism^0(\RR^n),$ $T_1=T,$ and $T_t$ having the canonical
relation $\Phi_t.$ Let $S_t$ be a family of parametrices.  Then, for $A \in
\Psitwo(\RR^n;o),$ we have
$$ \frac{d}{dt} T_t A S_t \cong [T_t'S_t, T_t A S_t] \bmod
\Psism^{-\infty}(\RR^n).
$$ Hence, setting $A(t) = T_t A S_t,$ we have
\begin{equation}\label{theode}
A'(t) \cong [P(t), A(t)] \bmod \Psism^{-\infty}(\RR^n),
\end{equation}
with 
$$
A(0) \in \Psitwo^{m,l}(\RR^n;0)
$$ having the desired wavefront and symbol properties, by the properties of
the calculus $\Psitwo.$  A priori, we have
$$
P(t) = T_t'S_t \in \Psism^{1}(\RR^n).
$$
However, since the canonical relation of $T_t$ is always the identity on
$o$, we can parametrize the FIOs $T_t$ by phase functions of the form
\begin{equation}\label{phases}
\phi(t,x,\theta) = (x-y)\cdot \theta+O(\theta^2);
\end{equation}
thus, $\pa_t \phi = O(\theta^2).$  
Differentiating \eqref{FIO},
we see that there are two terms in $T_t',$ coming from differentiation
of the phase $\phi$ and the amplitude $b$; the latter gives a term in
$T_t'S_t$ lying in $\Psism^{0}(\RR^n)\subset \Psitwo^{0,0}(\RR^n;o).$ The
former, by \eqref{phases}, has amplitude $h^{-1} O(\theta^2) b,$ i.e.\ is
of the form $h^{-1}$ times an order zero FIO with symbol vanishing to
second order on $(o\times o) \cap\diag.$ Consequently, the contribution to
$T_t'S_t$ from this term is an element of $\Psism^{1}(\RR^n)$ with
principal symbol vanishing on the diagonal to second order. By property
\eqref{hcalcintwocalc}, one order of vanishing yields
$$
P(t) \in \Psitwo^{1,0}(\RR^n;\lag);
$$ the second order of vanishing additionally gives\footnote{One can see
this simply by writing the total symbol in the form $\sum
h^{-1}\xi_i\xi_j+O(1)$ and lifting it to $[T^*\RR^n \times [0,1); o \times
0].$}
$$ \sigmatwo(P(t))=0 \text{ on } SN(\lag).
$$ Thus the ODE \eqref{theode} can be solved for $A(t)\in
\Psitwo^{m,l}(\RR^n;o)$ order-by-order, with a remainder in $\residual^l$
(which can be integrated away).  Moreover, the principal symbol on
$SN(\lag)$ is manifestly constant, as the Hamilton vector field of $P(t)$
vanishes there.  Thus 
$$
\sigmatwo (T A S)\restrictedto_{SN(\lag)} = \sigmatwo(A)\restrictedto_{SN(\lag)}.
$$ On the other hand, on $\pssp\backslash SN(\lag),$ the corresponding
statement follows from the usual semiclassical Egorov theorem and property
\eqref{awayfromlag} of the calculus.

The statement about microsupports is likewise straightforward from the ODE.
\end{proof}

\begin{definition}\label{def:two-micro-lag}
Suppose $X$ is a manifold without boundary and $\lag$ is a Lagrangian
submanifold of $T^*X$ such that the restriction of the bundle
projection, $\pi_{\lag}:\lag\to X$, is proper. We say that a family
of operators
$A=A_h:\CI(X)\to\CI(X),$ $h\in(0,1)$, with properly supported
Schwartz kernel $K_A\in\CmI(X\times X\times[0,1);\pi_R^*\Omega X)$, is in
$\Psitwo^{m,l}(X;\lag)$ if
\begin{enumerate}
\item\label{it1}
for $Q\in\Psism^0(X)$ with $\WF'(Q)\cap\lag=\emptyset$, $QA,AQ\in\Psism^m(X)$,
\item\label{it2}
for each point $q\in\lag$ and neighborhood $\U\subset T^*X$ of $q$
symplectomorphic, via a canonical transformation $\Phi$, to a neighborhood
of $q'\in o\subset T^*\RR^n$, mapping $\lag$ to $o$, and for each
semiclassical Fourier integral operator $T$ with canonical relation $\Phi$
elliptic in $\U$, with
parametrix $S$, and for each $Q,Q'\in\Psism^0(X)$ with $\WF'(Q),\WF'(Q')
\subset \U$, $TQAQ'S\in\Psitwo^{m,l}(\RR^n;o)$,
\item\label{it3}
if $Q,Q'\in\Psism^0(X)$ with $\WF'(Q)\cap\WF'(Q')=\emptyset$, $T$, resp $T'$,
elliptic semiclassical FIOs
mapping a neighborhood of $\WF'(Q)$, resp.\ $\WF'(Q')$, to a neighborhood
of $o$, and $\lag $ to $o$, with parametrices $S$, resp.\ $S'$, then
$TQAQ'S'\in\residual^l$.
\end{enumerate}
\end{definition}

\begin{definition}\label{def:quantize}(A global quantization map.)
Let $\{\U_j:\ j\in J\}$ be an open cover of $\lag$ such that
\begin{enumerate}
\item
$\bar \U_j\subset T^*X$ is compact for each $j$,
\item
for each
$K\subset X$ compact, $\pi^{-1}(K)\cap \U_j=\emptyset$ for all but finitely
many $j$, where $\pi:T^*X\to X$ is the bundle projection.
\item
for each $j$ there is a canonical transformation $\Phi_j$ from $\U_j$
to an open set $\U'_j$ in $T^*\RR^n$ mapping $\lag$ to $o$, with inverse
$\Psi_j$, and a
semiclassical Fourier integral operator $T_j$ elliptic in $\U_j$ with
canonical relation $\Phi_j$ with parametrix $S_j$.
\end{enumerate}
(Such an open cover exists because each point in $\lag$
has a neighborhood satisfying (1) and (3), and as $\pi_\lag$ is
proper, (2) can be fulfilled as well.)
Let $J_*=J\cup\{*\}$ be a disjoint union, and let $\U_*=T^*X\setminus \lag$,
so $\{\U_j:\ j\in J_*\}$ is an open cover of $T^*X$. Let
$\{\chi_j,\ j\in J_*\}$ be a
subordinate partition of unity. For $a\in S^{m,l}(\symbsp)$, let
\begin{equation}\label{eq:global-quant}
\Op(a)=\sum_{j\in J_*} S_j\Opw(\Psi_j^*(\chi_j a))T_j.
\end{equation}
(Here one could use $\Opl$ or $\Opr$ instead of $\Opw$ to obtain another
quantization.)
\end{definition}

\begin{definition}
With $(X;\lag)$, $A\in\Psitwo^{m,l}(X;\lag)$
as above, if $O\subset T^*X\setminus\lag$, $Q\in\Psism^0(X)$,
$\WF'(Q)\cap\lag=\emptyset$, $\WF'(\Id-Q)\cap O=\emptyset$, then
\begin{equation*}
\sigmatwo_{m,l}(A)|_O=\sigma_h(QA)|_O,\ \WF'(A)\cap O=\WF'(QA)\cap O,
\end{equation*}
while if $Q,Q',S,T$ as
in \eqref{it2} of Definition~\ref{def:two-micro-lag}, $O$ a neighborhood
of $q$ such that $O\cap\WF'(\Id-Q)=\emptyset$,
$O\cap\WF'(\Id-Q')=\emptyset$, then
\begin{equation*}\begin{split}
&\sigmatwo_{m,l}(A)|_O=\Phi^*
\big(\sigmatwo_{m,l}(TQAQ'S)|_{\Phi(O)}\big),\\
& \WF'(A)\cap O
=\Phi^{-1}(\WF'(TQAQ'S)\cap
\Phi(O)).
\end{split}\end{equation*}
\end{definition}

It is easy to check that in the overlap regions, the various cases give the
same classes of operators. For instance, if $Q,Q'\in\Psism^0(X)$ with
$\WF'(Q),\WF'(Q') \subset \U$ as in \eqref{it2}, but
$\WF'(Q)\cap\WF'(Q')=\emptyset$, then taking $G,G'\in\Psism^0(X)$ such that
$\WF'(G)\cap\WF'(G')=\emptyset$, $\WF'(\Id-G)\cap\WF'(TQS)=\emptyset$,
$\WF'(\Id-G')\cap\WF'(TQ'S)=\emptyset$,
$TQAQ'S-GTQAQ'SG'\in\Psism^{-\infty}(X)$ as
$(\Id-G)TQS\in\Psism^{-\infty}(X)$ by properties of the standard
semiclassical calculus, etc., so $TQAQ'S\in\Psitwo^{m,l}(\RR^n;o)$, and
Lemma~\ref{lemma:WFp} shows that $GTQAQ'SG'\in\residual^l$, so
$TQAQ'S\in\residual^l$ as well.  Thus, in the overlap, where both make
sense, the cases \eqref{it2} and \eqref{it3} are equivalent.

Correspondingly,
$\Op$ indeed maps into $\Psitwo(X;\lag)$:
\begin{equation*}
TQS_j\Opw(\Psi_j^*(\chi_j a))T_jQ' S\in\Psitwo^{m,l}(\RR^n;o),
\ TQS_j\Opw(\Psi_j^*(\chi_j a))T_j Q'S'\in\residual^l,
\end{equation*}
as $TQS_j$ and
$T_j Q'S$ are semiclassical Fourier integral operators
preserving the zero section. It is
similarly easy to check that the principal symbol and operator wave front
sets are well-defined (one only needs to check on $SN(\lag),$ as away from
this face they agree with the corresponding semiclassical quantities).

The proof of the properties (i)-(xi) follows from the case $\lag=o$
using the semiclassical FIOs as in the definition,
Proposition~\ref{prop:FIO-o} and Lemma~\ref{lemma:off-lag}.

If $\lag$ is a torus,
there is an improved quantization map (for symbols supported sufficiently
close to $\lag$) for which the full asymptotic
formula for composition is given by the formula from Weyl calculus.
First, suppose that $X=\TT^n$, and $\lag$ is the zero section. Let
$\{\phi_i:\ i\in I\}$ be a partition of unity subordinate to a finite
cover of $X$ by coordinate charts $(O_i,F_i)$, $G_i=F_i^{-1}$,
such that the transition maps $F_i\circ G_j$
between coordinate charts are all given by translations in $\RR^n$, and
let $\psi_i\equiv 1$ on a neighborhood of $\supp\phi_i$.
Then for $a\in S^{m,l}(\symbsp)$, define
$$
\Op(a)=\sum_i F_i^*\psi_i\Opw(\phi_i a)\psi_i G_i^*.
$$
For $b\in S^{m,l}(\symbsp)$
supported in $T^*_{O_i\cap O_j}X$,
$$
\Opw(b)-(F_j \circ G_i)^*
\Opw(b)(F_i\circ G_j)^*\in\residual^l.
$$
It follows that for $a\in S^{m,l}(\symbsp)$,
$$
\Op(a)^*-\Op(\bar a)\in\residual^l,
$$
{\em with the adjoint taken with respect to a translation invariant measure,}
and the composition formula in the Weyl calculus
(cf.\ \cite[Theorem 7.3 et seq.]{Dimassi-Sjostrand})
holds for the {\em global}
quantization map $\Op$: for $a\in S^{m,l}(\symbsp)$, $b\in S^{m',l'}(\symbsp)$,
\begin{equation}\begin{split} \label{weylcomposition}
\Op(a)\Op(b) = \Op\left( \sum
\frac{h^{\abs{\alpha+\beta}}(-1)^{\abs{\alpha}}}{(2i)^{\abs{\alpha+\beta}}\alpha!\beta!}\left((\pa_x^{\alpha}\pa_\xi^\beta
a) (\pa_\xi^{\alpha}\pa_x^\beta b)
\right)\right)+E
\end{split}\end{equation}
where the sum is a Borel sum, and $E \in
\Psitwo^{-\infty, l+l'}(\TT^n;o).$

Now, if $\lag$ is a Lagrangian torus in $T^*X$, there may not exist in
general a globally defined Fourier integral operator from a neighborhood of
$\lag$ in $T^*X$ to a neighborhood of the zero section in $T^*\TT^n$, even
though the underlying canonical relation $\Phi$ exists: such a choice may
in general be obstructed by both the Maslov bundle and Bohr-Sommerfeld
quantization conditions.  In fact, as we are conjugating, a multi-valued
FIO suffices, as noted by Hitrik and Sj\"ostrand \cite{Hitrik-Sjostrand1}
(see also \cite{Silva} in the non-semiclassical case). We can phrase this
slightly differently, with the notation of Definition~\ref{def:quantize},
locally identifying $\TT^n$ with $\RR^n$, by choosing an open cover of
$\lag$ by open sets $\U_j$ ($j\in J$) with $\U_k\cap \U_j$ contractible for
all $k,j\in J$, and choosing Fourier integral operators $T_j$ associated to
the canonical relation $\Phi|_{\U_j}$ mapping from the open subsets $\U_j$
of $T^*X$ to $T^*\TT^n$ (mapping $\lag$ to the zero section), such that
\begin{enumerate}
\item
$T_j$ is
unitary (modulo $O(h^\infty)$) on a neighborhood of $\supp\chi_j$, $j\in J$,
i.e.\ $T^*_jT_j-\Id\in\Psism^{-\infty}(X)$ microlocally near $\supp\chi_j$
(so we can take $S_j=T_j^*$),
\item
on $\supp\chi_k\cap \supp\chi_j$, $T^*_k T_j\in\Psism^0(X)$ is
given by multiplication by $c_{kj}e^{i\alpha_{kj}/h}$ for some constants
$c_{kj}\neq 0$ and $\alpha_{kj}$ (modulo $O(h^\infty)$) and for $j,k\in J$
(so not necessarily so if $j,k\in J_*$),
\item
and finally,
for each $j\in J$, in local coordinates $x$ on $\TT^n$ and $y$ on $X$, in
which the measure is $\abs{dx}$, resp.\ $\abs{dy}$,
$T_j$ is given by an oscillatory integral of the form
$$
T_j u(x)=h^{-(n+N)/2}\int e^{i\phi(x,y,\theta)/h}t_j(x,y,\theta,h)u(y)\,dy
\,d\theta
$$
where $\phi$ parameterizes the Lagrangian corresponding to $\Phi|_{\U_j}$,
and $t_j|_{C_\phi\times\{0\}}$ has constant argument,
where $C_\phi=\{(x,y,\theta):
d_\theta\phi(x,y,\theta)=0\}$, so the graph of $\Phi|_{\U_j}$ is
$\{((x,d_x\phi(x,y,\theta)),(y,-d_y\phi(x,y,\theta)):\ (x,y,\theta)\in
C_\phi\}$.
\end{enumerate}
Such a
choice of the $T_j$ exists (see
\cite[\S 2]{Hitrik-Sjostrand1}), and the last condition implies (see
\cite[Proposition~2.1]{Hitrik-Sjostrand1}) that for
$p\in h^{-m}\CI_c(T^*X\times
[0,1))$,
\begin{equation}\label{eq:weyl-conjugation-prop}
T_j^*\Opw((\Phi_{\U_j}^{-1})^*p) T_j-\Opw(p)\in \Psism^{m-2}(X).
\end{equation}
Then \eqref{eq:global-quant} gives a global
quantization map.

As $T_kS_j\in\Psism^0(\TT^n)$
is
given by multiplication by $c_{jk}^{-1}e^{-i\alpha_{jk}/h}$
(modulo $O(h^\infty)$), if $b\in S^{m,l}(\symbsp)$
is supported in $\Phi(\U_j)\cap\Phi(\U_k)$,
then $S_j\Opw(b)T_j-S_k\Opw(b)T_k\in\Psitwo^{-\infty,l}(X)$, for
$T_kS_j\Opw(b)T_jS_k-\Opw(b)\in\residual^l$.
Thus for $a\in S^{m,l}(\symbsp)$, $b\in S^{m',l'}(\symbsp)$
with
$\esssupp a,\esssupp b\subset T^*X\setminus\supp\chi_*$,
\begin{equation*}\begin{split}
\Op(a)\Op(b)&=\sum_{j,k\in J}T_j^*\Opw(\Psi^*(\chi_j a))T_jT_k^*\Opw(\Psi^*
(\chi_k b))T_k\\
&\cong \sum_{j,k\in J}T_j^*T_jT_k^*\Opw(\Psi^*(\chi_j a))\Opw(\Psi^*
(\chi_k b))T_k\\
&\cong \sum_{j,k\in J}T_k^*\Opw(\Psi^*(\chi_j a))\Opw(\Psi^*
(\chi_k b))T_k\ \text{modulo}\ \residual^{l+l'},
\end{split}\end{equation*}
using properties (1) and (2) of the $T_j$, so using the symplectomorphism
invariance of the Weyl composition formula,
\begin{equation}\begin{split} \label{weylcomposition-lag}
\Op(a)\Op(b) = \Op\left( \sum
\frac{h^{\abs{\alpha+\beta}}(-1)^{\abs{\alpha}}}{(2i)^{\abs{\alpha+\beta}}\alpha!\beta!}\left((\pa_x^{\alpha}\pa_\xi^\beta
a) (\pa_\xi^{\alpha}\pa_x^\beta b)
\right)\right)+E
\end{split}\end{equation}
where the sum is a Borel sum, computed in any local coordinates, and $E \in
\Psitwo^{-\infty, l+l'}(X;\lag).$
Since
modulo $\residual$ composition is microlocal, it suffices if one
of $a,b$ satisfies the support condition.

In addition, if $a$ satisfies the
support condition, then $\Op(a)^*-\Op(\bar a)\in\residual^l$, so replacing
$\Op$ by
$$
\Op'(a)=\frac{1}{2}(\Op(a)+\Op(a)^*)+\frac{1}{2}(\Op(a)-\Op(\bar a)),
$$
for real-valued $a$, $\Op'(a)$ is self-adjoint.
Thus, we have the following
result:

\begin{proposition}\label{prop:lag-weyl}
Suppose that $\lag$ is a Lagrangian torus in $T^*X$, with $\pi_\lag:\lag\to X$
proper.
Then there exists a neighborhood $\U$ of $\lag$ in $T^*X$ and
a quantization map $\Op:S^{m,l}(\symbsp)\to
\Psitwo^{m,l}(X;\lag)$ satisfying all properties listed in
Section~\ref{sec:calculus}, and such that, in addition,
\begin{enumerate}
\item
if $O$ is a coordinate chart in $X$ in which the volume form is given by
the Euclidean measure, then for
$p\in h^{-m}\CI_c(T^*O\times[0,1))$
with $\esssupp p\subset \U$,
\begin{equation}\label{eq:weyl-vs-improved}
\Op(p)-\Opw(p)\in \Psism^{m-2}(X).
\end{equation}
\item
for
$a\in S^{m,l}(\symbsp)$, $b\in S^{m',l'}(\symbsp)$ with either
$\esssupp a\subset \U$ or $\esssupp b\subset \U$,
\eqref{weylcomposition-lag} holds,
\item
if $a$ is real-valued with
$\esssupp a\subset\U$, then $\Op(a)$ is
self-adjoint,
\item
for any $a$ satisfying
$\esssupp a\subset\U$, $\Op(a)^*-\Op(\bar a)\in\residual^l$.
\end{enumerate}
\end{proposition}

\section{Real principal type propagation}

Recall that we let $\pssp$ denote the space $[T^*X; \lag]$ on which
principal symbols of operators live and $\symbsp=[T^*X \times [0,1);
    \lag\times 0]$ the space for total symbols.
If $P\in \Psitwo^{m,l}(X;\lag)$ has real principal symbol $p,$ its Hamilton
vector field $\bar{\sH}$ is a vector field on $\pssp,$ conormal to $SN(\lag)
=\pa\pssp.$  If
$\sbdf$ is a boundary defining function for $SN(\lag)\subset \pssp,$ then
the appropriately rescaled Hamilton vector field
$$
\sH = \sbdf^{l-m+1} \bar{\sH}
$$ is a smooth vector field on $\pssp$, tangent to its boundary,
$SN(\lag).$ In particular, if a point in an orbit of $\sH$ is in
$\pa\pssp$, the whole orbit is in $\pa\pssp$.

The following result is the corresponding real principal
type propagation theorem.

\renewcommand{\theenumi}{\Roman{enumi}}
\begin{theorem}\label{theorem:0}
Let $P\in \Psitwo^{m,l}(X;\lag)$ have real principal symbol $p.$ If
$u\in I^{-\infty}_{(r)}(\lag)$, $P u=f$ then
\begin{enumerate}
\item 
$\displaystyle{\WFtwo^{k,r}(u)\backslash\WFtwo^{k-m,r-l}(f)
\subset \Sigmatwo(P) \equiv \{\sigmatwo(P)=0\}.}$
\item
$\displaystyle{\WFtwo^{k,r}(u)\backslash
\WFtwo^{k-m+1,r-l}(f)}$ is invariant under the
Hamilton flow of $p$ inside $\Sigmatwo(P).$
\end{enumerate}
\end{theorem}
\renewcommand{\theenumi}{\arabic{enumi}}

\begin{proof}
This is just the usual real principal type propagation away from
the boundary of $\pssp$, so we only need to consider points at
$\pa\pssp$.

The proof of the first part follows from the existence of elliptic
parametrices (property \eqref{ellipticparametrix} of the calculus).

The proof of the flow invariance follows the outline of H\"ormander's
classic commutator proof of the propagation of singularities for operators
of real principal type \cite{MR48:9458}, hence we give only a sketch here.
(See also \cite{MR95k:58168} for an account of essentially the same proof
in the setting of a different pseudodifferential calculus.)

Pick $q \in \Sigmatwo(P) \subset SN(\lag).$ Let $\sH$ denote the
Hamilton vector field of $P,$ and assume that
\begin{equation}\label{blah}
\exp(r_0 \sH) q \notin \WFtwo^{\alpha,r} (u)
\end{equation}
and
\begin{equation}\label{blahblah}
\exp(t \sH) q \notin \WFtwo^{\alpha-1/2,r} (u),\ \exp(t \sH) q \notin \WFtwo^{\alpha-m+1,r-l} (f)\quad \text{ for } t \in [0, r_0];
\end{equation}
we will show that for $r_0>0$ sufficiently small (depending on $\sH$, but
not depending on $u$) \eqref{blah}--\eqref{blahblah} imply that
\begin{equation}\label{conclusion}
\exp(t \sH) q \notin \WFtwo^{\alpha,r} (u)\quad \text{ for } t \in [0, r_0].
\end{equation} We obtain the corresponding result with the interval $[0,r_0]$ in
\eqref{blah} replaced by $[-r_0, 0]$ by applying the result to the operator
$-P.$ If $\alpha<k,$ we can then iterate this argument to obtain the
desired result.

To prove
that \eqref{blah}--\eqref{blahblah} imply \eqref{conclusion},
let $\chi_0(s)=0$ for $s\leq 0$, $\chi_0(s)=e^{-M/s}$ for $s>0$ ($M>0$ to be
fixed),
let $\chi\geq 0$ be
a smooth non-decreasing
function on $\RR$ with $\chi=0$ on $(-\infty, 0]$ and $\chi=1$ on
$[1,\infty),$ with $\chi^{1/2}$ and $(\chi')^{1/2}$ both smooth.
Let $\phi$ be a cutoff function supported in
$(-1,1).$ As $\sH$ is a vector field tangent to the boundary
of $\pssp$, assuming that $\sH$ does not vanish at $q\in
\pa\pssp$ (otherwise there is nothing to prove), we can
choose local coordinates $\rho_1,\dots, \rho_{2n}$ on $[T^*X; \lag]$
centered at $q$ in which $\sH=\pa_{\rho_1}$, and the boundary is defined
by $\rho_{2n}=0$, so we may take $\sbdf=\rho_{2n}$. We choose $r_0$ so that
$\exp(t\sH)q$, $t\in[-r_0,2r_0]$, remains in a compact subset of the coordinate
chart given by the $\rho_j$.
Let $\rho' = (\rho_2,\dots,
\rho_{2n}),$ and set
\begin{equation*}\begin{split}
a=
&\sbdf^{-r+l/2+\alpha-(m-1)/2}
\phi^2(\lambda^2|\rho'|^2) \chi_0(\lambda\rho_1+1)
\chi(\lambda(r_0-\rho_1)-1)\\
&\in
\sbdf^{-r+l/2+\alpha-(m-1)/2}\CI(\pssp).
\end{split}\end{equation*}
Then $a$ has support in the region where we have assumed regularity,
provided $\lambda$ is chosen large enough since
$|\rho'|<\lambda^{-1}$,
$\rho_1\geq-\lambda^{-1}$, $\rho_1\leq r_0-\lambda^{-1}$
on $\supp a$.  Let
$A\in \Psitwo^{\alpha-(m-1)/2,r-l/2}(X;\lag)$
have principal symbol $a$;
recall that the principal symbol of an element of
$\Psitwo^{\alpha-(m-1)/2,r-l/2}(X;\lag)$ arises by factoring out
$h^{-(\alpha-(m-1)/2)}$
from the `amplitude' $a_{\tot}$ (with $A=\Opl(a_{\tot})$
locally), explaining the power
of $\sbdf$ appearing above.

Noting that the weight function $\sbdf=\rho_{2n}$ has vanishing derivative
along $\sH$, we compute
\begin{equation*}\begin{split}
\sigmatwo_{2\alpha,2r}(-i(A^*AP-P^*A^*A))&=\sigmatwo_{2\alpha,2r}(-i[A^*A,P])
+\sigmatwo_{2\alpha,2r}(-iA^*A(P-P^*))\\
&=\sH a^2-i\sigmatwo_{m-1,l}(P-P^*) a^2=b^2-e^2+a^2 q,
\end{split}\end{equation*}
where $b,e\in \sbdf^{r-\alpha}\CI(\pssp)$
(arising from symbolic terms in
which $\sH$ has been applied to $\chi_0(\lambda\rho_1+\lambda^{-1})^2,$
and $\chi(\lambda(r_0-\rho_1))^2$ respectively), $q\in\sbdf^{l-m+1}\CI(\pssp)$
(arising from $P-P^*$).
Thus, in view of the principal symbol short exact
sequence,
$$
-i(A^*AP-P^*A^*A) = B^*B-E^*E+A^*QA+R
$$ with $B,E \in \Psitwo^{\alpha,r}(X;\lag)$
$Q\in\Psitwo^{m-1,l}(X;\lag)$, having $\WF'(B)$, etc., given by
$\esssupp b$, etc., $\sigmatwo_{\alpha,r}(B)=b$, etc., and $R \in
\Psitwo^{2\alpha-1,2r}(X;\lag),$ with $\WF'(R)\subset\esssupp a$.
We thus find that $\WF'(E)$ is contained in $|\rho'|<\lambda^{-1}$,
$\rho_1\in[r_0-2\lambda^{-1},r_0-\lambda^{-1}]$, hence (for sufficiently
large $\lambda$) in the
complement of $\WFtwo^{\alpha,r} (u),$ so $\norm{E u}$ is uniformly bounded,
as $\WF'(R)$ is contained in
$|\rho'|<\lambda^{-1}$,
$\rho_1\in[-\lambda^{-1},r_0-\lambda^{-1}]$, so $|\ang{R u,u}|$ is also
uniformly bounded,
and $B$ is elliptic on $\exp(t \sH)
q,$ for $t \in [0, r_0-2\lambda^{-1}],$
with
$a\leq CM\lambda^{-1}b$.  Thus
$$
\norm{B u}^2 \lesssim \abs{\ang{R u, u}}
+\abs{\ang{A u, Af}}+\abs{\ang{A u,Q^* A u}}+ \norm{E u}^2.
$$ Let $T\in \Psitwo^{(m-1)/2,l/2}(X;\lag)$ be elliptic on a neighborhood
of $\WF'(A)$ and let $T'$ be a
parametrix; we rewrite $\abs{\ang{A u, Af}} \lesssim \delta \norm{T A u}^2+
\frac{1}{4\delta}\norm{T' Af}^2$ modulo residual errors;
$ TA\in\Psitwo^{\alpha,r}(X,\lag)$ with principal
symbol
a smooth non-vanishing multiple of $a$, so for $\delta$
sufficiently small the first term may be absorbed into $\norm{B u}^2$ (as
$\sqrt{b^2-c^2 a^2}$ is smooth for small $c>0$) modulo a residual term,
while the second is uniformly bounded as $h \to 0$ by our
assumption \eqref{blahblah}. On the other hand,
$\abs{\ang{ TA u,  T'QAu}} \lesssim \norm{ T A u}^2+
\norm{ T'QAu}^2$ modulo residual errors, and
$ TA, T'QA\in\Psitwo^{\alpha,r}(X,\lag)$ with principal
symbol
a smooth non-vanishing multiple of $a$.
For $M$ sufficiently small, both terms can be absorbed into $\norm{Bu}^2$
(modulo a term that can be absorbed into $R$).
\end{proof}

\section{Propagation of $\WFtwo$ on Invariant Tori in Integrable Systems}\label{section:propagation}

Let $P\in \Psism^0(X).$ Assume that
\begin{equation}\label{Phypoth1}P \text{ has real principal symbol $p$}
\end{equation}
with Hamilton vector field denoted $\sH,$
and assume that 
\begin{equation}\begin{split}\label{Phypoth2}
\sH &\text{ is completely integrable in a
neighborhood of}\\
&\text{a compact Lagrangian invariant torus }\lag\subset\{p=0\}.
\end{split}\end{equation}
Let $(I_1,\dots,I_n,\theta_1,\dots,\theta_n)$ be the associated
action-angle variables.  Without loss of generality we may translate the
action coordinates so that $\lag$ is defined by $I_i=0$ for all
$i=1,\dots,n.$ Let $\omega_i = \partial p/\partial I_i$ and $\omega_{ij} =
\partial^2 p/\partial I_i \partial I_j.$ Let $\ombar_i$ and $\ombar_{ij}$
denote the corresponding quantities restricted to $\lag$ (where they are
constant).  We introduce coordinates on $[T^*X; \lag]$ by setting
$$
\rho= \abs{I} = \bigl(\sum I_j^2\bigr)^{\frac 12},\quad \hI_j= \frac{I_j}{\abs{I}}.
$$
The front face of the blown-up space is defined by $\rho=0$ and is
canonically identified with the spherical normal bundle $SN(\lag).$

The real principal symbol assumption on $P$ means that (with respect to the
inner product on $L^2(X)$ given by {\em any} smooth density on $X$)
$P^*-P\in \Psism^{-1}(X)$, i.e.\ $P$ is self-adjoint to leading order. It
turns out that for our improved result we need at the very least that
$\sigma_{h,-1} (P^*-P)$ vanishes at $\lag$ with respect to some inner
product; unlike the statement that $P^*-P\in \Psism^{-1}(X)$, this depends
on the choice of an inner product. So we assume from now on that $X$ has a
fixed density $\nu$ on it (e.g.\ a Riemannian density).  This density $\nu$
in turn yields a trivialization of the bundle of half-densities on
$X.$ As
is well-known, this yields a canonically defined \emph{subprincipal symbol}
$\sub A$ for a semiclassical pseudodifferential operator $A\in\Psism^m(X).$
In our (semiclassical) setting, $\sub P$ can be defined using
the Weyl quantization---cf.\ the improved symbol invariance statement in
Lemma~\ref{lemma:diffeo-inv} (see also, for instance, \cite{Hormander:vol3} in
the non-semiclassical case). To obtain the subprincipal symbol we thus choose a
coordinate system in which the Euclidean volume form agrees with the fixed
one (this can always be arranged by changing one of the coordinates, while
fixing the others); writing $A=\Opw (a)$ in these coordinates, we have
\begin{equation}\label{eq:subpr-symbol-expansion}
a=\sigma_h(A)+h \sub_h(A)+O(h^2).
\end{equation}
As $\Opw(a)^*=\Opw(\bar a)$ (adjoint taken with respect to the Euclidean
volume form), if $\sigma_{h,m}(A)$ is real,
$$
\sigma_{h,m-1}(A-A^*)=2i\Im\sub_h(A),
$$
so $A-A^*\in\Psism^{m-2}(X)$ if and only if $\sub_h(A)$ is real.

We now impose a weakened self-adjointness condition on $P$, namely that
$P-P^*\in\Psism^{-2}(X)$ with respect to the fixed density, i.e.\
$\sub_h(P)$ is real; we further assume that $\sub_h(P)$ is constant on
$\lag$:
\begin{equation}\label{Phypoth3}
\sub_h(P) \text{ is real on }T^*X,\text{ and it is constant on } \lag.
\end{equation}
In fact, the slightly weaker assumption
\begin{equation}\label{Phypoth3p}
\sub_h(P) \text{ is real and constant on } \lag
\end{equation}
would suffice; one would need to take care of $P-P^*$ much as in the
proof of Theorem~\ref{theorem:0}: we assume \eqref{Phypoth3} as it covers the
cases of interest.

We recall that as $\lag$ is characteristic, we have $P \in
\Psitwo^{0,-1}(X;\lag),$ and the principal symbol $\sigmatwo_{0,-1}(P)$ near
$SN(\lag)$ is
\begin{equation}\label{princsymb}
\sum\ombar_j I_j+ \frac 12 \sum \ombar_{ij} I_i
I_j+O(I^3) =\rho \sum\ombar_j \hI_j+ \frac{\rho^2}{2}\sum \ombar_{ij} \hI_i
\hI_j+O(\rho^3) 
\end{equation}
The Hamilton vector field is thus
$$
\sH = \sum \ombar_j \pa_{\theta_j}+ \rho \sum \ombar_{ij} \hI_i
\pa_{\theta_j}+\rho^2 \sH'
$$ with $\sH'$ smooth on $\pssp=[T^*X; \lag]$ and tangent to the boundary,
$SN(\lag).$

\begin{theorem}\label{theorem:1}
Assume that $u \in L^2$ and $Pu\in O(h^\infty)L^2,$ where $P$ and $\lag$
satisfy \eqref{Phypoth1},\eqref{Phypoth2},\eqref{Phypoth3}.  Then for each
$k$ and each $l\leq 0,$ $\WFtwo^{k,l} (u)\cap SN(\lag)$ is invariant under
\begin{equation}
\sH_1 = \sum \ombar_j \pa_{\theta_j}.
\end{equation} Also, for each $k,$ and each $l\leq -1/2,$ $\WFtwo^{k,l} (u) \cap
SN(\lag)$ is invariant under
\begin{equation}\label{H2}
\sH_2 = \sum \ombar_{ij} \hI_i \pa_{\theta_j}
\end{equation}
\end{theorem}
\begin{remark}
In fact, for the first part of the theorem it suffices to adopt the weaker
hypotheses that
$$
\WFtwo^{k+1,l+1} (Pu) = \emptyset
$$
and for the second part
$$
\WFtwo^{k+1,l+2} (Pu) = \emptyset;
$$
for instance, it certainly suffices to require $Pu \in h^s L^2$ where
$s \geq \max(k+1,l+1)$ or $s \geq \max(k+1,l+2)$ rspectively.
\end{remark}

The rest of this section will be devoted to a proof.

Invariance under $\sH_1$ follows from Theorem~\ref{theorem:0}, as $h^{-1} P\in
\Psitwo^{1,0}(X;\lag),$ and $\sH_1$ is the restriction of the Hamilton
vector field to $SN(\lag);$ note that $\WF^{-\infty,l} (u) = \emptyset$ for
each $l\leq 0$ as $u \in L^2.$

Invariance under $\sH_2$ requires considerable further discussion.

Given $\zeta \in SN(\lag),$ suppose that we know that $\WF^{k-1/2,l} u$ is
invariant and that $\zeta \notin \WFtwo^{k,l} (u).$ We need to show that the
$\sH_2$ orbit through $\zeta$ is disjoint from $\WFtwo^{k,l} (u);$ the general
case can be obtained from this argument by the usual iteration.  We know that
the closure of the $\sH_1$-orbit of $\zeta$ is a torus $T_\zeta\subset
SN(\lag),$ and that $T_\zeta\cap \WFtwo^{k,l} u=\emptyset$ by
$\sH_1$-invariance.

We extend the vector field $\sH_1$ to a neighborhood of $\pa\pssp$ using
the coordinates $(I,\theta)$ as above.  Thus, the closure of each orbit of
$\sH_1$ near $\pa\pssp$ is still a torus. Let
$\sH_2$ be defined on
a neighborhood of $SN(\lag)$ by
$$
\sH_2=\rho^{-1}(\sH-\sH_1);
$$ this naturally agrees with \eqref{H2} on $SN(\lag).$ Then
$[\sH_1,\sH_2]=0$ near $\lag$ as the principal symbol of $P$
is independent of $\theta$ there, so $\sH$, $\sH_1$ are linear combination
of the $\pa_{\theta_j}$ with coefficients depending on $I$ only, and
$\rho=|I|$.
Note that $I$ is constant
along the flow of $\sH_1$ and
$\sH_2$, so for $\delta>0$ small, the $(\sH_1,\sH_2)$ joint flow from
$\rho<\delta$ stay in this prescribed
neighborhood of $SN(\lag)$ for all times, on which $I$ and $\theta$ are thus
defined.

Now let $a_0\in \CI(\pssp)$ be supported in the complement of $\WFtwo u,$
but sufficiently close to $SN(\lag)$ (so that $I$, etc., are defined
on a neighborhood of the $(\sH_1,\sH_2)$-flowout of $\supp a_0$),
with $a_0$ having smooth square root, and
$$
\sH_1 a_0=0.
$$
For $\delta\in(0,1)$ small, let
$$
a_1 = -\int_0^\delta (2-s) a_0 \circ \exp(-s \sH_2) \, ds
$$
so that
$$
\rho \sH_2 (a_1) = \rho\int_0^\delta a_0 \circ
\exp(-s \sH_2) \, ds +(2-\delta)\rho a_0\circ
\exp(- \delta\sH_2)-2\rho a_0 = \int_0^\delta b_s^2\, ds+c^2-d^2
$$ with $b^2,c^2,d^2$ real and vanishing to first order at $\rho=0,$
given by the three terms in the middle expression in the displayed formula.
Note also
that by construction
$$
\sH_1 (a_1) = 0
$$
since $\sH_1$ and $\sH_2$ commute.

Now let
\begin{equation}\label{eq:a-tot-def}
a=h^{-2k}|I|^{-(2l+1)+2k}a_1,
\end{equation}
and let $\Op$ be the quantization given in Proposition~\ref{prop:lag-weyl} corresponding
to the local symplectomorphism $\Phi=(\theta,I)$ near $\lag$, mapping
a neighborhood of $\lag$ to a neighborhood of the zero section of $\TT^n.$  Thus
$$
A=\Op(a)\in \Psitwo^{2k,2l+1} (X)
$$
is selfadjoint, with
$$
\sigmatwo(A)=h^{-2k}|I|^{-(2l+1)+2k}a_1,
$$
and if $P=\Op(p)$ microlocally near $\U$, $p\in \CI_c(T^*X\times[0,1))$, then
\begin{equation}\begin{split}\label{weylcommutator}
i&((h^{-1}P)^*A-A(h^{-1}P)\cong ih^{-1}(\Op(\bar p)\Op(a)-\Op(a)\Op(p))\\
&\cong
\Op\left( \sum
\frac{h^{\abs{\alpha+\beta}-1}(-1)^{\abs{\alpha}}}{(2i)^{\abs{\alpha+\beta}}\alpha!\beta!}\left(
(\pa_\theta^{\alpha}\pa_I^\beta\bar p) (\pa_I^{\alpha}\pa_\theta^\beta a)
-(\pa_\theta^{\alpha}\pa_I^\beta a) (\pa_I^{\alpha}\pa_\theta^\beta p)
\right)\right)\\
&\qquad\qquad\qquad \text{modulo}\ \residual^{2l+1}.
\end{split}\end{equation}
Here we used the symplectomorphism invariance of the Weyl
composition formula in order to utilize the action-angle coordinates
(which simply undoes the already used symplectomorphism invariance
in obtaining \eqref{weylcomposition-lag}).

Note that by \eqref{eq:weyl-vs-improved} and
\eqref{eq:subpr-symbol-expansion},
if $P=\Op(p),$ $p=p_0+hp'_1+O(h^2)$ then
\begin{equation}\label{betteregorov}
\begin{aligned}
\sigma_h(P) &= p_0\\
\sub_h(P) &= p'_1.
\end{aligned}
\end{equation}
(The main novelty here is the formula for the subprincipal symbol.)

Let $B_s,C,D \in \Psitwo^{k,l}(X;\lag)$ have symbols $b_s,c,d$
respectively.  We may easily compute in the usual manner (since the weights
commute with $P$ to leading order):
$$
\sigmatwo(i((h^{-1}P)^*A-A(h^{-1}P))) = \rho\sH_2(a_1)+2a\Im\sub_h(P)
=\rho\sH_2(a_1)
$$
by \eqref{Phypoth3}.
Thus we have
\begin{equation}\label{ourcommutator}
i((h^{-1}P)^*A-A(h^{-1}P)) = \int_0^\delta B_s^* B_s \, ds+C^*C-D^*D+R
\end{equation} with $R \in \Psitwo^{2k-1,2l+1}(X;\lag),$ and $\WFtwo'(R) \subset
\WFtwo'(A).$ Note, then, that a priori $R$ has higher order than $C^*C$ in the
second index, as the invariance of $\sigmatwo(A)$ along the $\sH_1$ flow
yields a vanishing of the \emph{principal} symbol of the ``commutator'' at
$SN(\lag),$ but not necessarily of lower-order terms.
However, the
use of the special quantization $\Op$ gives us a better result:

\begin{lemma}
We may decompose
$$
R=R_1+R_2 \in
\Psitwo^{2k-1,2l}(X;\lag)+\Psitwo^{-\infty, 2l+1}(X;\lag).
$$
\end{lemma}

\begin{proof}[Proof of Lemma]
The
subprincipal symbol of $P$ is a real constant on $\lag;$
let $\mu$ denote
this constant.  Thus we have $P=\Op(p),$ $A=\Op(a),$ with
\begin{equation}\label{p}\begin{split}
p&=p_0
+\mu h
+O(I h)+O(h^2)\\
& =\sum\ombar_j I_j+ \frac{1}{2}\sum \ombar_{ij} I_i I_j
+\mu h+O(\rho_{\ff}^2 \rho_{\sidef})\equiv p_0+p_1,\\
p_0&=\sum\ombar_j I_j+ \frac 12 \sum \ombar_{ij} I_i
I_j+O(I^3),
\end{split}\end{equation}
where $p_0$ is independent of $h$, and $a$ as in \eqref{eq:a-tot-def}.
Here $O(\rho_{\sidef}^{-r}\rho_{\ff}^{-s})$ stands for an element of
$S^{r,s}(\symbsp),$ and as usual $\rho_{\ff}$ denotes a boundary defining
function for the front face, and $\rho_{\sidef}$ a boundary defining
function for the side face of the blowup $\symbsp=[T^*X;\lag].$

Now we use \eqref{weylcommutator}. Writing
$p=p_0+p_1$ (as in \eqref{p}), the terms arising by replacing
$p$ by $p_0$ in \eqref{weylcommutator} with $\alpha=\beta=0$ cancel, while
the terms with $|\alpha+\beta|=1$ give $\Op(H_{p_0}a)\in\Psitwo^{2k,2l}(X)$,
which differs from $\int_0^\delta B_s^* B_s \, ds+C^*C-D^*D$ by an element
$\tilde R$
of $\Psitwo^{2k-1,2l}(X)$ (as they are both in $\Psitwo^{2k,2l}(X)$ and
have the same principal symbol). Thus,
$R$ is obtained by taking the terms in \eqref{weylcommutator}
arising from $p_1$, as well as those arising from $p_0$ with
$\abs{\alpha+\beta}>1,$ along with the remainder of the formula and $\tilde R$.

We now examine \eqref{weylcommutator} in coordinates on
$\symbsp=[T^*X;\lag]$ given locally by $H=h/I_1,$ $I_1$,
$\tilde I=I'/I_1=(I_2/I_1,\dots,I_n/I_1),$ and $\theta_1,\dots,\theta_n$
(these are of course valid only in one part of the corner of the blowup,
but other patches are obtained symmetrically).  Thus, $H$ is a defining
function for $\sidef$ and $I_1$ for $\ff.$ By the analogous computation
to \eqref{vflifts}, all terms have, a priori, the same conormal order at
$SN(\lag)$ (the terms all have asymptotics $I_1^{-2l-1}$ in the
coordinates employed in \eqref{vflifts}, with powers of $H$ ascending from
$H^{2k+1}$). However, since $p$ is actually a \emph{smooth} function
on $T^*X\times[0,1),$ hence we have $\pa_I^\gamma p=O(1)$ for
all $\gamma$ (with the above notation, so further $\pa_\theta$ derivatives have
the same property), so if $|\beta|\geq 1$,
$$
h^{|\alpha+\beta|-1}(\pa_\theta^{\alpha}\pa_I^\beta
p) (\pa_I^{\alpha}\pa_\theta^\beta a)\in
S^{2k+1-|\alpha|-|\beta|,2l+2-|\beta|}(\symbsp).
$$
Moreover by \eqref{p},
$$
\pa_{\theta}^\gamma p = O(\rho_{\sidef} \rho_{\ff}^2)
=O(H I_1^2)\quad \text{ if }\abs{\gamma}>0.
$$
Hence for $|\beta|\leq|\alpha|$, $|\alpha|\geq 1$,
$$
h^{|\alpha+\beta|-1}(\pa_\theta^{\alpha}\pa_I^\beta
p) (\pa_I^{\alpha}\pa_\theta^\beta a)\in
S^{2k-|\alpha|-|\beta|,2l}(\symbsp).
$$
We thus conclude that the terms of the
form
$$
h^{|\alpha+\beta|-1}(\pa_\theta^{\alpha}\pa_I^\beta
p) (\pa_I^{\alpha}\pa_\theta^\beta a)
$$ with $|\alpha+\beta|>1$ (which thus have either $|\beta|\geq 2$, or
$|\alpha|\geq 1$, $|\beta|\leq|\alpha|$) in the sum are in fact
$O(H^{-2k+1} I_1^{-2l}).$ An analogous calculation holds if we interchange
the role of $\alpha$ and $\beta$ (as well as if we complex conjugate), and
we conclude that modulo the terms with $|\alpha+\beta|=1$, we can Borel sum
the right hand side of \eqref{weylcommutator} to a symbol in
$$
S^{-2k+1,-2l}(\symbsp).
$$
\emph{This concludes the proof of the Lemma.}
\noqed \end{proof}

The remainder of the proof of the theorem is as follows.

Pairing \eqref{ourcommutator} with $u$ we obtain for all $h>0,$
$$ -2 \Im\ang{Au, h^{-1}Pu} = \int_0^\delta \norm{B_su}^2\, ds
+\norm{Cu}^2-\norm{Du}^2+\ang{R_1 u,u}+ \ang{R_2 u,u}.
$$ The assumption of absence of $\WFtwo^{k,l}(u)$ at $\zeta$ (and hence its
$\sH_1$-orbit) controls $\ang{Du,u}$ uniformly as $h\downarrow 0.$ The
assumption of absence of $\WFtwo^{k-1/2,l}(u)$ on the whole microsupport of
$A$ controls $\ang{R_1 u,u}.$ The assumption $u \in L^2$ controls $\ang{R_2
u,u}$ since $R_2 \in \Psitwo^{-\infty,2l+1} = \residual^{2l+1} \subset
\residual^0$ since $l \leq -1/2.$ Thus, since the left side is
$O(h^\infty),$ we obtain absence of $\WFtwo^{k,l}(u)$ on the elliptic set
of $C$, hence on the time-$\delta$ flowout of $\sH_2$ (for any small
$\delta$). Hence we obtain $\sH_2$-invariance of $\WFtwo^k(u).$ (A
corresponding argument works along the backward flowout.)\qed

\section{Consequences for Spreading of Lagrangian Regularity}

Recall that an invariant torus in an integrable system (with the notation
of \S\ref{section:propagation}) is said to be \emph{isoenergetically
  nondegenerate} if
\begin{equation}\label{bigmatrix}
\Omega=
\begin{pmatrix}
\omega_{11}& \dots & \omega_{1n} & \omega_1 \\
\vdots & \ddots & \vdots & \vdots\\
\omega_{n1} & \dots & \omega_{nn} & \omega_n\\
\omega_{1} & \dots & \omega_{n} & 0\\
\end{pmatrix}
\end{equation}
is a nondegenerate matrix.  We recall from \cite{Wunsch:Integrable} that a
somewhat trivial example of a system in which the invariant tori are
nondegenerate is when $P=h^2\Lap-1$ on $S^1\times S^1;$ we may take
$\lag$ to be, for instance, $\{\xi_1=1,\ \xi_2=1\}.$ Here $\Lap$ is the
nonnegative Laplacian, and $\xi_i$ are the fiber variables dual to $x_i$ in
$T^*(S^1\times S^1).$ A considerably less trivial example is the spherical
pendulum, where all tori are isoenergetically nondegenerate except for
those given by a codimension-one family of exceptional energies and angular
momenta---see Horozov \cite{Horozov1,Horozov2} for the proof of the
nondegeneracy and the description of the exceptional tori.

\begin{definition}
A distribution $u$ is Lagrangian on a closed set $F \subset \lag$ if there
exists $A \in \Psism(X),$ elliptic on $F,$ such that $Au$ is a Lagrangian
distribution with respect to $\lag.$
\end{definition}

We recall that in \cite{Wunsch:Integrable}, it was shown that under the
hypotheses of Theorem~\ref{theorem:1},\footnote{In \cite{Wunsch:Integrable}
the subprincipal symbol assumption was in fact stronger: it was assumed to
vanish.} and if $\lag$ is assumed to be isoenergetically nondegenerate,
then local Lagrangian regularity on $\lag$ is invariant under the Hamilton
flow of $P$ on $\lag,$ and, additionally, Lagrangian regularity on a small
tube of closed bicharacteristics implies regularity along the
bicharacteristics inside it.

We now prove the following, generalizing the results of
\cite{Wunsch:Integrable}.

\begin{corollary}\label{cor1}
Assume that the hypotheses of Theorem~\ref{theorem:1} hold and that,
additionally, $\lag$ is isoenergetically nondegenerate.  If $u$ is
Lagrangian microlocally near any point in $\lag$ relative to $L^2$ then $u$
is globally Lagrangian with respect to $\lag$ in a microlocal neighborhood
of $\lag,$ relative to $h^\epsilon L^2$
for all $\ep>0.$
\end{corollary}
\begin{remark}
If the initial data for the wave equation on $\RR^n$ is smooth in an
annulus, the solution is smooth near the origin at certain later times.
This remark is to H\"ormander's propagation of singularities theorem as
Corollary~\ref{cor1} and the results of \cite{Wunsch:Integrable} are to
Theorem~\ref{theorem:1}: in both settings the crude statements about
singular supports are deducible from a much finer microlocal theorem.
\end{remark}
\begin{remark}
An example from \cite{Wunsch:Integrable} shows that the hypothesis of
isoenergetic nondegeneracy is necessary: without it, there do exist
quasimodes that are Lagrangian only on parts of $\lag.$

The reader may wonder if nowhere Lagrangian quasimodes are in fact
possible, given the hypotheses of the theorem.  An example is as follows:
consider
$$
P = h^2 \Lap-1
$$
on $S^1\times S^1$ (with $\Lap$ the nonnegative Laplacian).  Consider the
sequence
$$
u_k = e^{i(k^2x_1+kx_2)},\ k \in \NN.
$$
taking the sequence of values $h=h_k=k^{-1}(1+k^2)^{-1/2}$ gives
$$
P_{h_k} u_k=0.
$$ Now the $u_k$'s are easily verified (say, by local semiclassical Fourier
transform) to have semiclassical wavefront set in the Lagrangian
$\lag=\{(x_1,x_2,\xi_1=1,\xi_2=0)\}.$ On the other hand, the operator $h
D_{x_2}\in \Psi_h(S^1\times S^1)$ is characteristic on $\lag$ and we have,
for the sequence $h=h_k,$
$$
(h D_{x_2})^m u_k = (h k)^m u_k = (1+k^2)^{-m/2} u_k.
$$ Thus, $u$ certainly does not have iterated regularity under the
application of $h^{-1} (h D_{x_2}),$ hence is not a semiclassical
Lagrangian distribution.  (We recall, though, from \cite{Wunsch:Integrable}
that the hypotheses of the Corollary are satisfied in this case.)
\end{remark}

\begin{proof}
Let $x \in \lag,$ and assume that $u$ enjoys Lagrangian regularity at $x,$
relative to $L^2,$
i.e.\ $(x,\xi) \notin \WFtwo^{\infty,0}(u)$ for all $\xi \in SN_x(\lag).$  By
Theorem~\ref{theorem:1}, $\WFtwo^{\infty,-1/2}(u)$ is invariant under the flow
$$
\sum \ombar_j \pa_{\theta_j}
$$
and
$$
\sum \ombar_{ij} \xi_i \pa_{\theta_j}.
$$
and hence under any linear combination of them.  In other words, at given
$\vec{\xi} = (\xi_1,\dots,\xi_n)^t,$ $\WFtwo(u)$ is invariant under the
flow along the vector field
$$
V(\xi, s) = (\vec{\xi}^t, s)\cdot \Omega \cdot
\begin{pmatrix}\pa_{\theta}\\ 0\end{pmatrix}
$$ for all $s \in \RR,$ where $\Omega$ is given by \eqref{bigmatrix}.  By
isoenergetic nondegeneracy, the set of $\xi$ such that $V(\xi, s)$ has
rationally independent components for some $s \in \RR$ has full measure,
i.e.\ in particular such values of $\xi$ are dense in $SN_x(\lag).$ As the
closure of the orbit along such a vector field is all of $\lag\times \{\xi\},$ and as
$\WFtwo$ is a closed set, we find that for a dense set of $\xi \in
SN_x(\lag),$ $\WFtwo(u) \cap (\lag \times \{\xi\}) = \emptyset.$\footnote{Our
notation reflects the fact that $N(\lag)$ is a trivial bundle.}  Again by
closedness of $\WFtwo(u),$ we now find that $\WFtwo^{\infty,-1/2}(u) =
\emptyset.$  Since $u \in L^2,$ an interpolation yields
$$
\WFtwo^{\infty,-\ep}(u) = \emptyset.\qed
$$
\noqed\end{proof}

\bibliography{all}
\bibliographystyle{amsplain}
\end{document}


\maketitle

Nalini Anantharaman has pointed out to the authors that the proof of Corollary 6.2 of \cite{twomicro} is incomplete, as the assertion of the penultimate sentence is wrong.  We have been unable to complete the proof of the corollary as stated, but note that the incomplete proof does yield several interesting results, including a new proof of the somewhat mysterious results of \cite{Wunsch:Integrable}.  It is our hope that adopting some of the second microlocal methods used by Anantharaman-Maci\'a \cite{Anantharaman-Macia1} in the special setting of flat tori may lead to stronger results.

We emphasize that the gap in the proof of Corollary 6.2, which is an
application of the main propagation results of the paper, does not
affect any of the other results of \cite{twomicro}. In
particular the key results, Theorem 4.1 and Theorem 5.1, are valid with
unchanged proofs.  Also, as we show below, Corollary 6.2 does hold in dimension two.

In order to state our general replacement for Corollary 6.2, we first introduce a definition.  (We continue with the notation of \cite{twomicro}.)
\begin{definition}
Let $p \in \lag,$ and $V \subset T_p\lag \cong \RR^n$ a vector subspace.  Let $\overline{\exp_p(V)}$ denote the closure of the flowout starting at $p$ along all constant vector fields in $V.$
\end{definition}

Our result now describes how local Lagrangian regularity propagates,
based on our results on propagation of the second microlocal wavefront
set.  The idea is that regularity on a set $U\subset \lag$ spreads
automatically to certain points $p\in \lag:$ By Theorem~5.1 of
\cite{twomicro}, for any point $q\in\WFtwo(u)$ lying over a point $p$,
$\WFtwo(u)$ would automatically fill out the orbit closure under
$\Span\{\sH_1,\sH_2\}$ (with orbits emanating from $q$), hence if if all of these orbit closures
intersect $U,$ there cannot be any $\WFtwo$ over $p.$ Although
$\sH_1(q)$ regarded as a vector in (i.e.\ pushed forward to) $T_p\lag$
is independent of $q$, in the same sense $\sH_2(q)$ varies with $q$,
requiring a hypothesis on general two-dimensional $V$ below:

\begin{proposition}
Assume that the hypotheses of Theorem 5.1 hold and that, additionally, $\lag$ is isoenergetically nondegenerate.
Let $U \subset \lag$ and let 
\begin{equation*}\begin{split}
\widetilde{U}=\big\{p \in \lag:\ \overline{\exp_p(V)}\cap U \neq
\emptyset
\text{ for every two-dimensional } &V \subset \RR^n\\
&\text{ containing } \sum \ombar_j \pa_{\theta_j} \big\}.
\end{split}\end{equation*}
If $u$ is
locally Lagrangian on $U \subset \lag$ relative to $L^2$ then $u$
is locally Lagrangian on $\widetilde{U}$ relative to $h^{-\epsilon} L^2$ for all $\ep>0.$
\end{proposition}
\begin{proof}
As stated in the proof of Corollary 6.2,
$\WFtwo^{\infty,-1/2}(u)$ is invariant under the flow along
$$
\sH_1=\sum \ombar_j \pa_{\theta_j}
$$
and
$$
\sH_2= \sum \ombar_{ij} \hI_i \pa_{\theta_j}
$$
and hence under any linear combination of these vector fields.
We also remark that $\WFtwo^{\infty,-1/2}(u)$ is contained in the characteristic set
$$
\Sigma_2 = \{\sum \ombar_j \hI_j=0\}.
$$
We claim that $\sH_1, \sH_2$ are never linearly dependent on $\Sigma_2:$ if they were linearly dependent at $\hI = \xi,$ we would have by rescaling $\xi$ appropriately:
$$
\sum \ombar_{ij} \xi_i =\ombar_j
$$
for each $j,$ hence
$$
\begin{pmatrix}
\ombar_{11}& \dots & \ombar_{1n} & \ombar_1 \\
\vdots & \ddots & \vdots & \vdots\\
\ombar_{n1} & \dots & \ombar_{nn} & \ombar_n\\
\ombar_{1} & \dots & \ombar_{n} & 0\\
\end{pmatrix}
  \begin{pmatrix}
\xi_1\\ \vdots \\ \xi_n \\ -1
\end{pmatrix}=  \begin{pmatrix}
0\\ \vdots \\ 0 \\ \sum\xi_j \ombar_j
\end{pmatrix}.
$$
On $\Sigma_2$, this vector would vanish, contradicting isoenergetic nondegeneracy.  Thus, $\Span\{\sH_1,\sH_2\}$ is always two-dimensional.

If $p \in \widetilde{U},$ then for every value of $\xi$, we have $(p, \xi)\notin \WFtwo^{\infty, -1/2} u,$ since if $(p,\xi)$ did lie in the wavefront set, then by closedness, all of $\overline{\exp_p(V)}$ would lie in the wavefront set with $V =\Span\{\sH_1,\sH_2\},$ and this set intersects $U$ by hypothesis.  Thus, $u$ enjoys Lagrangian regularity with respect to $h^{-1/2}L^2$ on $\widetilde{U}.$  Since $u \in L^2,$ an interpolation yields the desired regularity.
\end{proof}

Consequently, we do recover the result stated in \cite{twomicro} in the special case $n=2:$
\begin{corollary}
Let $n=2.$   If $u$ is
locally Lagrangian relative to $L^2$ near some point $p \in \lag,$ then $u$
is locally Lagrangian on relative to $h^{-\epsilon} L^2$ on all of $\lag.$
\end{corollary}

We also recover the results of \cite{Wunsch:Integrable}, showing that Lagrangian regularity on an annular neighborhood of a closed orbit ``fills in'' the interior:
\begin{corollary}
Assume that all the $\ombar_j$ are rationally related, and let $\gamma$ denote a closed orbit under $\sH_1.$  Let $N_\delta$ denote a delta-neighborhood of $\gamma$ with respect to the flat metric on $\lag.$
If $u$ is locally Lagrangian relative to $L^2$ on all of $N_{2\delta} \backslash N_\delta$ then  $u$
is locally Lagrangian relative to $h^{-\epsilon} L^2$ on $\gamma.$
\end{corollary}
\begin{proof}
We simply choose $U=N_{2\delta} \backslash N_\delta$ and remark that for $p \in \gamma,$ any two-dimension orbit-closure through $p$ passes through $U$.
\end{proof}

\bibliography{all}
\bibliographystyle{amsplain}

%% file: symbolspace.customized.eepic
\setlength{\unitlength}{0.00029167in}
\begingroup\makeatletter\ifx\SetFigFont\undefined%
\gdef\SetFigFont#1#2#3#4#5{%
  \reset@font\fontsize{#1}{#2pt}%
  \fontfamily{#3}\fontseries{#4}\fontshape{#5}%
  \selectfont}%
\fi\endgroup%
{\renewcommand{\dashlinestretch}{30}
\begin{picture}(5049,6939)(0,-10)
\put(2274.028,3464.476){\arc{1505.144}{4.5965}{7.8700}}
\put(2224,3462){\ellipse{824}{1500}}
\path(2262,4212)(2262,6912)
\path(2262,2712)(2262,312)(2262,12)
\path(2787,3987)(5037,6387)
\path(1887,3087)(12,1437)
\dottedline{15}(2594,3794)(2787,3987)
\put(3162,3462){ff}
\put(2637,5787){sf}
\end{picture}
}